\documentclass[12pt]{amsart}
\usepackage{amsmath,amsfonts,amsthm,amsopn,amssymb,mathtools,stmaryrd,mathrsfs}
\usepackage{cite,marginnote}
\usepackage{pdfsync}

\usepackage{color,enumitem,graphicx}
\usepackage[colorlinks=true,urlcolor=blue,
citecolor=red,linkcolor=blue,linktocpage,pdfpagelabels,
bookmarksnumbered,bookmarksopen]{hyperref}
\usepackage[english]{babel}

\usepackage[left=2.75cm,right=2.75cm,top=3cm,bottom=3cm]{geometry}

\usepackage[hyperpageref]{backref}
\usepackage[]{xcolor}

\def\sideremark#1{\ifvmode\leavevmode\fi\vadjust{\vbox to0pt{\vss
 \hbox to 0pt{\hskip\hsize\hskip1em
 \vbox{\hsize2.1cm\tiny\raggedright\pretolerance10000
  \noindent #1\hfill}\hss}\vbox to15pt{\vfil}\vss}}}%

\numberwithin{equation}{section}

\makeindex
\newtheorem{theorem}{Theorem}[section]
\newtheorem{proposition}[theorem]{Proposition}
\newtheorem{lemma}[theorem]{Lemma}
\newtheorem{remark}[theorem]{Remark}
\newtheorem{example}[theorem]{Example}
\newtheorem{corollary}[theorem]{Corollary}
\newtheorem{definition}[theorem]{Definition}

\newcommand{\R}{\mathbb R}

\newcommand{\bt}{\begin{theorem}}
\newcommand{\et}{\end{theorem}}
\newcommand{\bl}{\begin{lemma}}
\newcommand{\el}{\end{lemma}}
\newcommand{\bd}{\begin{definition}}
\newcommand{\ed}{\end{definition}}
\newcommand{\bc}{\begin{corollary}}
\newcommand{\ec}{\end{corollary}}
\newcommand{\bp}{\begin{proof}}
\newcommand{\ep}{\end{proof}}
\newcommand{\bx}{\begin{example}}
\newcommand{\ex}{\end{example}}
\newcommand{\bi}{\begin{exercise}}
\newcommand{\ei}{\end{exercise}}
\newcommand{\bo}{\begin{proposition}}
\newcommand{\eo}{\end{proposition}}
\newcommand{\br}{\begin{remark}}
\newcommand{\er}{\end{remark}}
\newcommand{\be}{\begin{equation}}
\newcommand{\ee}{\end{equation}}
\newcommand{\ba}{\begin{align}}
\newcommand{\ea}{\end{align}}
\newcommand{\bn}{\begin{enumerate}}
\newcommand{\en}{\end{enumerate}}
\newcommand{\bg}{\begin{align*}}
\newcommand{\bcs}{\begin{cases}}
\newcommand{\ecs}{\end{cases}}

\newcommand{\bean}{\begin{eqnarray*}}
\newcommand{\eean}{\end{eqnarray*}}

\renewcommand{\leq}{\leqslant}
\renewcommand{\geq}{\geqslant}

\title[From $H^1$-critical to super-critical quasilinear Schr\"{o}dinger equations]{Blow-up phenomena and asymptotic profiles passing from $H^1$-critical to super-critical quasilinear Schr\"{o}dinger equations}

\author[D.~Cassani]{Daniele Cassani$^\text{1}$}
\author[Y.~Wang]{Youjun Wang$^\text{2}$}

\address[D. Cassani]{\newline\indent Dip. di Scienza e Alta Tecnologia
\newline\indent
Universit\`{a} degli Studi dell'Insubria
\newline\indent and
\newline\indent RISM--Riemann International School of Mathematics
\newline\indent Villa Toeplitz, Via G.B. Vico, 46 -- 21100 Varese}
\email{\href{mailto:Daniele.Cassani@uninsubria.it}{Daniele.Cassani@uninsubria.it}}

\address[Y.~Wang]{\newline\indent Department of Mathematics
\newline\indent
South China University of Technology
\newline\indent Guangzhou 510640, P. R. China
\newline\indent and
\newline\indent Dip. di Scienza e Alta Tecnologia
\newline\indent
Universit\`{a} degli Studi dell'Insubria
\newline\indent
via Valleggio 11, 22100 Como,Italy}
\email{\href{mailto:scyjwang@scut.edu.cn}{scyjwang@scut.edu.cn}}

\thanks{(1) Corresponding author: \texttt{Daniele.Cassani@uninsubria.it}}
\thanks{(2) The second named author is supported by the Fundamental Research Funds for the Central Universities (No.~2018MS59) and Natural Science Foundation of Guangdong (No.~2018A0303130196)}

\date{\today}

\begin{document}

\maketitle

\begin{abstract}
We study the asymptotic profile, as $\hbar\rightarrow 0$, of positive solutions to
$$-\hbar^2\Delta
u+V(x)u-\hbar^{2+\gamma}u\Delta u^2=K(x)|u|^{p-2}u,\ \ x\in
\mathbb{R}^N
$$ where $\gamma\geq 0$ is a parameter with relevant physical interpretations, $V$ and $K$ are given potentials and $N\geq 5$. We investigate the concentrating behavior of solutions when $\gamma>0$ and, differently form the case $\gamma=0$ where the leading potential is $V$, the concentration is here localized by the source potential $K$. Moreover, surprisingly for $\gamma>0$ we find a different concentration behavior of solutions in the case $p=\frac{2N}{N-2}$ and when $\frac{2N}{N-2}<p<\frac{4N}{N-2}$. This phenomenon does not occur when $\gamma=0$.
\end{abstract}

\section{Introduction}\label{introduction}

\noindent We are concerned with blow-up phenomena for positive solutions to the following class of quasilinear Schr\"{o}dinger
equations
\begin{equation}\label{maineq} -\hbar^2\Delta
u+V(x)u-\hbar^{2+\gamma}u\Delta u^2=K(x)|u|^{p-2}u,\ \ x\in
\mathbb{R}^N, N\geq 5
 \end{equation} where $\hbar>0$ is the adimensionalized Planck constant,  $\gamma\in\R$ is a parameter which is relevant in several applications in Physics  for which we refer to \cite{Kur81,KWR01}, and which we assume here to be positive, $V$ and $K$ are given potentials, for the moment real continuous functions,
and the nonlinearity in the range $\frac{2N}{N-2}\leq p<\frac{4N}{N-2}$. We refer to \cite{Sre90,Br-Er-01,Br-Er-03,NA77} and
references therein  for the Physics context of (\ref{maineq}).

\noindent The existence of nontrivial solutions, in particular ground states for (\ref{maineq}) has been intensively studied in recent years throughout a very extensive literature, among which let us mention \cite{Liu03, Liu04, SV01}. Though it is not possible to give exhaustive references on the subject, let us recall a few results which are strictly related to our problem.

\noindent For
   semi-classical states of (\ref{maineq}), namely $\gamma=0$ and $\hbar\rightarrow 0$,
assume $2<p<\frac{4N}{N-2}$, $N\geq 3$,  $K(x)\equiv 1$,
 and that $V:\mathbb{R}^N\rightarrow \mathbb{R}$  is H\"{o}lder continuous and
satisfying the following conditions:
$0<V_0<\inf\limits_{x\in \mathbb{R}^N}V(x)$ and there is a  bounded  open set $\Lambda$ such that
$0<a:=\inf\limits_{x\in \Lambda}V(x)<\min\limits_{x\in \partial\Lambda}V(x)$. Then,  the existence of localized solutions concentrating near $\Omega:=\{x\in \Lambda:
V(x)=a\}$ has been obtained in \cite{DAN09,Gloss01} and, by scaling properties, as $\hbar\rightarrow 0$ the limit equation turns out to be the following quasilinear autonomous Schr\"odinger equation
\begin{equation}\label{e00}-\Delta
u+au- u\Delta u^2=|u|^{p-2}u,\ \ x\in
\mathbb{R}^N.
 \end{equation}
We refer to  \cite{do10, Wang-Zou,He13,do09, Yang13} for related results.

\noindent Notice that the scaling invariance of \eqref{maineq} breaks down as soon as $\gamma>0$. Recently in \cite{CWZ20}, it has been proved that in this context both the cases $\gamma=0$ and $\gamma>0$ have similar concentration behavior. However, the limit equation for $\gamma>0$ is different form the case $\gamma=0$ and turns out to be the following semilinear Schr\"odinger equation
 \begin{equation}\label{e11}-\Delta
u+au=|u|^{p-2}u,\ \ x\in
\mathbb{R}^N.
 \end{equation}

\noindent As we are going to see, this fact will play a crucial role in studying the blow-up profile of solutions to \eqref{maineq}. Indeed, loosely speaking one expects solutions can be localized along suitable normalized truncations and translations of ground states to the limit equation (\ref{e00}) or (\ref{e11}). Here the situation is completely different from the case $K\equiv 1$ and $\gamma=0$, as a proper normalized, translated and rescaled solution will concentrate around critical points of the potential $K$.

\noindent It is well known from \cite{AFM05,WZ97, BS08} that for the nonautonomous semilinear Schr\"{o}dinger equation
\begin{equation*}\label{maineq3} -\hbar^2\Delta
u+V(x)u=K(x)|u|^{p-2}u,\ \ x\in
\mathbb{R}^N,
 \end{equation*}
the function
$$\mathcal{A}(x)=[V(x)]^{\frac{p+2}{p}-\frac{N}{2}}[K(x)]^{-\frac{2}{p}}$$ retains important
information for the concentrating behavior of solutions.  Remarkably, for our problem (\ref{maineq}) the external Schr\"odinger potential $V$ does not play any role in the blow-up phenomenon which is governed by the source potential $K$.

\noindent Another interesting phenomenon addressed in this paper is the different concentrating behavior which occurs passing form critical to supercritical nonlinearities in \eqref{maineq}. This is due to the fact that the limit equation, as $\hbar\rightarrow 0$, for (\ref{maineq}) changes passing from $\frac{2N}{N-2}< p<\frac{4N}{N-2}$ to $p=\frac{2N}{N-2}$. Surprisingly, in the critical case the limit equation turns out to be the zero mass semilinear Schr\"odinger equation. This fact to the best of our knowledge has not been observed before.

\noindent In order to state our main results, set
 $$v(x)=\hbar^{\frac{\gamma}{2}}u(\hbar^{1+\frac{(p-2)\gamma}{4}}x),$$ then equation (\ref{maineq}) turns into the following
\begin{equation}\label{maineq1} - \Delta v+\hbar^{\frac{(p-2)\gamma}{2}}V(\hbar^{1+\frac{(p-2)\gamma}{4}}y)v- v\Delta v^2= K(\hbar^{1+\frac{(p-2)\gamma}{4}}y)|v|^{p-2}v,\ \ x\in
\mathbb{R}^N.
 \end{equation}

\noindent For simplicity set $\kappa=\hbar^{\frac{(p-2)\gamma}{2}}$, $\varepsilon=\hbar^{1+\frac{(p-2)\gamma}{4}}$ and denote $V(\hbar^{1+\frac{(p-2)\gamma}{4}}y),$ $ K(\hbar^{1+\frac{(p-2)\gamma}{4}}y)$ by $V_\varepsilon (y)$,  $K_\varepsilon(y)$, respectively.
Thus, equation (\ref{maineq1}) can  be  written in the following form
\begin{equation}\label{maineq2} - \Delta v+\kappa V_\varepsilon(y)v- v\Delta v^2= K_\varepsilon(y)|v|^{p-2}v,\ \ x\in
\mathbb{R}^N\ .
 \end{equation}

\noindent We assume the potential $V$ and $K$ satisfying the following conditions:
\begin{itemize}
\item [$(V)$]  $V\in C(\mathbb{R}^N,\mathbb{R})$, $0<\inf V(x)\leq V(x)\leq \sup V(x)<+\infty$;
\item [$(K)$]  $K:\mathbb{R}^N\rightarrow \mathbb{R}$ is H\"{o}lder continuous,  $0<\sup\limits_{x\in \mathbb{R}^N}K(x)<K_0$ and there is a  bounded  open set $\mathcal{O}$ such that
$$\max\limits_{x\in \partial \mathcal{O}}K(x)<m:=\sup\limits_{x\in \mathcal{O}}K(x)\ .$$
Set $\mathcal{M}:=\{x\in \mathcal{O}\ : \ K(x)=m\}$.
\end{itemize}

\medskip

\noindent Our main results are the following:
\begin{theorem} \label{Th1} Let $\gamma>0$, assume $(V)$, $(K)$ and
$\frac{2N}{N-2}\leq p<\frac{4N}{N-2},$ $N\geq 5$. Then, for sufficiently small $\varepsilon>0$,
 there exists a positive solution $v_\varepsilon$ of (\ref{maineq2}).
 \end{theorem}

\noindent The solution $v_\varepsilon$ obtained in Theorem \ref{Th1} is actually uniformly bounded with respect to $\varepsilon$.
As a consequence we will obtain  the blow-up profile of  solutions to the original equation \eqref{maineq}.

\noindent In Section 2, we prove some preliminary results, in particular we deal with the zero mass case and prove that the equation
 \begin{equation}\label{limit-eq}-\Delta u-u\Delta u^2=mu^p \end{equation}
has a unique positive radial solution $U$ which belongs to $D^{1,2}(\mathbb{R}^N)$.
Similarly  to \cite[Proposition 6.1]{AD14},    $v_\varepsilon\rightarrow U$ in $D^{1,2}(\mathbb{R}^N)\cap C_{loc}^2(\mathbb{R}^N)$, as $\varepsilon\rightarrow 0$. That is,
$$\hbar^{\frac{\gamma}{2}}u_\hbar(\hbar^{1+\frac{(p-2)\gamma}{2}}(\cdot-x_{\hbar}))\rightarrow   U(\cdot),$$
in $D^{1,2}(\mathbb{R}^N)\cap C_{loc}^2(\mathbb{R}^N)$, as $\hbar\rightarrow 0$.

\noindent Blow-up phenomena for  the autonomous version of problem (\ref{maineq}) (namely, $V(x)=\lambda>0$ and $K(x)\equiv1$) have been studied in \cite{AD14},
where in order to get the asymptotic profile of  the solution,
uniform estimates of the rescaled  ground state and energy estimates were established. However, their method can not be applied to deal with the non-autonomous problem (\ref{maineq}). In \cite{CW19}, the Lyapunov-Schmidt reduction method has been used to deal with the problem
\begin{equation}\label{maineq5} -\Delta
u+\varepsilon V(x)u- u\Delta u^2= u^p,\ \ u>0,\ \ \lim\limits_{|x|\rightarrow+\infty} u(x)=0\ \ x\in
\mathbb{R}^N.
 \end{equation}
 Assuming $V>0$, $V\in L^\infty$, $V(x)=o(|x|^{-2})$,  as $|x|\rightarrow +\infty$, the authors proved that for $\varepsilon$ sufficiently small,  problem (\ref{maineq5}) has a positive fast decaying solution provided $\frac{2N}{N-2}<p<\frac{4N}{N-2}$, $N\geq 3$.

\noindent Surprisingly, the limit equation for (\ref{maineq}) changes again when $p=\frac{2N}{N-2}$. Precisely,
let  $$v(x)=\hbar^{\frac{\alpha}{2}}u(\hbar^{1+\frac{(p-2)\alpha}{4}}x),\ \  \text{ for any } \ \ 0<\alpha<\gamma,$$    $\lambda=\hbar^{\frac{(p-2)\alpha}{2}}$, $\zeta= \hbar^{\gamma-\alpha}$, $\epsilon=\hbar^{1+\frac{(p-2)\alpha}{4}}$ and denote $V(\hbar^{1+\frac{(p-2)\alpha}{4}}y),$ $ K(\hbar^{1+\frac{(p-2)\alpha}{4}}y)$ by $V_\epsilon (y)$,  $K_\epsilon(y)$, respectively, then equation (\ref{maineq}) turns into the following equation
\begin{equation}\label{maineq4} - \Delta v+\lambda V_\epsilon(y)v- \zeta v\Delta v^2= K_\epsilon(y)|v|^{p-2}v,\ \ x\in
\mathbb{R}^N.
 \end{equation}
Note that $\lambda$, $\zeta\rightarrow 0$ as $\hbar\rightarrow 0.$

\noindent The solution to (\ref{maineq4}) is closely related to the (unique)
solution of the following zero mass mean field  limit equation \cite{BELI}
\begin{equation}\label{limit2} - \Delta v = mv^{\frac{N+2}{N-2}},\ \ x\in
\mathbb{R}^N, \ \ v>0, \ \ v(0)=\max v(x).
 \end{equation}

\noindent It is well known since \cite{TA} that equation (\ref{limit2}) possesses an explicit one parameter family of solutions given by
 $$U=(N(N-2)m)^{\frac{N-2}{4}}\left(\frac{\mu}{1+\mu^2|x|^2}\right)^{\frac{N-2}{2}},\ \ \mu>0.$$

\begin{theorem} \label{Th2} Assume $\gamma>0$, that $(V)$, $(K)$ hold and that
$\frac{2N}{N-2}<p<\frac{4N}{N-2},$ $N\geq 5$. Then, for sufficiently small $\hbar>0$, there exists a local maximum point $x_\hbar$ of $u_\hbar$ such that
$\lim\limits_{\hbar\rightarrow 0}\mbox{dist}(x_\hbar,\mathcal{M})=0$ and there exists a positive solution $u_\hbar$ of (\ref{maineq}) satisfying
$$u_\hbar(\cdot)= \hbar^{-\frac{\gamma}{2}}U(\hbar^{-1-\frac{(p-2)\gamma}{2}}\cdot-x_{\hbar})+\omega_\hbar(\cdot),$$ where $\omega_\hbar(\cdot)\rightarrow 0$ in $D^{1,2}(\mathbb{R}^N)\cap C_{loc}^2(\mathbb{R}^N)$  as $\hbar\rightarrow 0$, and $U$ is the unique   fast decay, positive and radial (least energy) solution of (\ref{limit-eq}). Moreover, when $p=\frac{2N}{N-2}$, under the above hypotheses we have
$$u_\hbar(\cdot)= \hbar^{-\frac{\alpha}{2}}\left[\frac{N(N-2)\mu\sqrt{m}}{1+\mu^2|\hbar^{-1-\frac{(p-2)\alpha}{2}}\cdot-x_{\hbar}|^2}\right]^{\frac{N-2}{2}}  +\omega_\hbar(\cdot),\quad \text{ for all }\quad 0<\alpha<\gamma\ .$$
\end{theorem}

\begin{remark} In Theorem \ref{Th2}, $\alpha=\gamma$ is not allowed. Indeed, there exist no fast decaying solutions to    (\ref{limit-eq}) if $p=\frac{2N}{N-2}$, as established in Theorem 1.1 of \cite{CW19}.
  \end{remark}

\medskip

\noindent Throughout this paper,  $C$ will denote a positive constant whose exact
value may change from line to line without affecting the overall result.

\section{Preliminaries}

\noindent In this Section we collect a few results, which we will use in the sequel, on the following zero mass equation
\begin{equation}\label{masszero} -
 \Delta
u-u\Delta u^2=m|u|^{p-2}u,\ \ \
x\in\mathbb{R}^N,  \end{equation}
where $2^*:=\frac{2N}{N-2}<p<2(2^*):=\frac{4N}{N-2},$ $N\geq 5$.

\noindent Uniqueness and non-degeneracy of positive solutions oto  (\ref{masszero}) have been completely solved in \cite{AD16}, see also \cite{CW19}. For convenience of the reader we recall below a few results we need in the sequel.

\noindent The energy functional related to equation $(\ref{masszero})$  is given by
$$I(u)=\frac{1}{2}\int_{\mathbb{R}^N}(1+2u^2)|\nabla u|^2dx-\frac{m}{p}\int_{\mathbb{R}^N}|u|^pdx$$
and it is well defined in the set $$ E=\left\{u\in D^{1,2}({\mathbb{R}^N}): \int_{\mathbb{R}^N} u^2|\nabla u|^2dx<+\infty \right\}.$$

\begin{theorem}[Theorem 1.1 in \cite{AD16} or Theorem 1.1 in \cite{CW19}]
Equation (\ref{masszero})  has a unique positive radial solution which belongs to $D^{1,2}(\mathbb{R}^N)$. In particular,  the ground state
 of (\ref{masszero}) is unique up to translations.
    \end{theorem}

\begin{lemma}[Lemma 2.1, \cite{YWA13}] \label{G}
Let $g(s)=\sqrt{1+2s^2}$ and $G(t)=\int_0^tg(s)ds$. Then, $G(t)$ is an odd smooth function as well as the inverse function $G^{-1}(t)$. Moreover, the following properties hold:
\begin{itemize}
\item [$(i)$]  $\displaystyle\lim\limits_{t\rightarrow 0}\frac{G^{-1}(t)}{t}=1$;\\

\item [$(ii)$]  $\displaystyle\lim\limits_{t\rightarrow +\infty}\frac{G^{-1}(t)}{\sqrt{t}}=\sqrt[4]{2}$;\\

\item [$(iii)$]   $   |G^{-1}(t)|\leq  |t|$  for all $t\in \mathbb{R}$;\\

\item [$(iv)$] $|G^{-1}(t)|^2$ is convex in $t$;\\

\item [$(v)$] $   |G^{-1}(t)|\leq   \sqrt[4]{2}\sqrt{|t|}$  for all $t\in \mathbb{R}$.
\end{itemize}
\end{lemma}

\noindent Next consider the following semilinear elliptic equation, which  in some sense is the dual problem of (\ref{masszero}):
\begin{equation}\label{mass1} -
 \Delta
v=m\frac{|G^{-1}(v)|^{p-2}G^{-1}(v)}{g(G^{-1}(v))},\ \ \
x\in\mathbb{R}^N,  \end{equation}

\noindent The energy functional corresponding to   $(\ref{mass1})$  is defined by
$$L_m(v)=\frac{1}{2}\int_{\mathbb{R}^N}|\nabla v|^2dx-\frac{m}{p}\int_{\mathbb{R}^N}|G^{-1}(v)|^pdx,$$
which is well defined in $D^{1,2}(\mathbb{R}^N)$ by Lemma \ref{G}, moreover $L_m(v)\in C^1$.

\noindent Solutions  to $v$ of (\ref{mass1}) satisfy the following Pohozaev identity
\begin{eqnarray}\label{Poh1}\frac{N-2}{2N}\int_{\mathbb{R}^N}|\nabla v|^2dx-\frac{m}{p}\int_{\mathbb{R}^N}|G^{-1}(v)|^pdx=0.\end{eqnarray}

\noindent Moreover, the ground state has  a mountain pass characterization, namely
$$L_m(U)=C_m=\inf_{\eta\in \Phi}\max_{t\in [0,1]}L_m(\eta(t)),$$where $\Phi=\{\eta\in C([0,1],D^{1,2}(\mathbb{R}^N)): \eta(0)=0, L_m(\eta(1))<0\}$, see Proposition 4.3 in \cite{AD14}.

\begin{theorem}[Propositions 2.6 and 3.2 in \cite{AD16}]
The following properties hold:
\begin{itemize}
\item [$(i)$]  (\ref{mass1}) has a unique fast decay positive radial solution $v(r)$, namely
$$\lim\limits_{r\rightarrow +\infty}r^{N-2}v(r)=c\in (0,+\infty);$$
\item [$(ii)$]Let $v\in D^{1,2}(\mathbb{R}^N)\cap C^2(\mathbb{R}^N) $  be a positive radially decreasing solution of (\ref{mass1}). Then there exists  $C>0$  such that
$$CA(r)(1-O(r^{-2}))\leq v(r)\leq CA(r), \ \ CA'(r)\leq v'(r)\leq CA'(r)(1-O(r^{-2})), $$
for sufficiently large $r$. Here
$$A(r)=\frac{1}{(N-2)|S^{N-1}|r^{N-2}}$$
 is the fundamental solution of $-\Delta$ on $\mathbb{R}^N$. In particular we have that
\begin{eqnarray}\label{estimate} \lim_{r\rightarrow +\infty}r^{N-2}v(r)=\frac{C}{(N-2)|S^{N-1}|},\ \ \lim_{r\rightarrow +\infty}r^{N-1}v'(r)=-\frac{C}{|S^{N-1}|}.\end{eqnarray}
\end{itemize}
    \end{theorem}

\begin{theorem} [Lemma 2.4 in \cite{AD16}] Suppose that $v=\int_0^ug(s)ds$. Then,
\begin{itemize}
\item [$(i)$]  $u\in E\cap C^2(\mathbb{R}^N)\Leftrightarrow v\in D^{1,2}(\mathbb{R}^N)\cap C^2(\mathbb{R}^N)$ ;
\item [$(ii)$]  $u$ is a positive solution of (\ref{masszero}) $\Leftrightarrow$ $v$ is a positive solution of (\ref{mass1}).
\end{itemize}
    \end{theorem}

\section{Proof of Theorems \ref{Th1} and \ref{Th2}}

\noindent We next consider  the following   quasilinear
Schr\"{o}dinger equation
\begin{equation}\label{main}
-\text{div}(g^2(u)\nabla u)+ g(u)g'(u)|\nabla u|^2+\kappa V_\varepsilon(x)u=K_\varepsilon(x)|u|^{p-2}u,\ \ x\in
\mathbb{R}^N,
\end{equation}
where $g(s)=\sqrt{1+2s^2}$. Direct calculations show that (\ref{main}) is equivalent to (\ref{maineq2}).

\noindent The energy
functional corresponding to  (\ref{main}) is given by
\begin{equation}\label{2.2}J_{\varepsilon}(u)=\frac{1}{2}\int_{\mathbb{R}^N}[g^2(u)|\nabla u|^2+\kappa V_\varepsilon(x)u^2]dx-\frac{1}{p}\int_{\mathbb{R}^N}K_\varepsilon(x)|u|^pdx.
\end{equation}
Note that $J_{\varepsilon}$ is not even well defined in $H^1(\mathbb{R}^N)$. However, it is well known since \cite{COJE,Liu03} that a suitable dual approach, hidden in change of variables, turns the energy functional to be smooth and well defined in a proper function space setting, see \cite{Shen13}, and also \cite{DAN09} for an Orlicz space approach. Here, the change of variables  $u=G^{-1}(v)$, yields the following smooth energy
\begin{eqnarray}\label{s2.2}P_{\varepsilon}(v)=\frac{1}{2}\int_{\mathbb{R}^N}\left(|\nabla v|^2+\kappa V_\varepsilon(x)|G^{-1}(v)|^2\right)dx-\frac{1}{p}\int_{\mathbb{R}^N}K_\varepsilon(x)|G^{-1}(v)|^pdx.
\end{eqnarray}
By Lemma \ref{G}, it is standard to check that $P_{\varepsilon}
\in C^1(H^1(\mathbb{R}^N),\mathbb{R})$.

\noindent The Euler-Lagrange equation associated to  $P_{\varepsilon}$ is the following
\begin{eqnarray}\label{S2.00} -
\Delta v+\kappa V_\varepsilon(x)\frac{G^{-1}(v)}{g(G^{-1}(v))}
=K_\varepsilon(x)\frac{|G^{-1}(v)|^{p-2}G^{-1}(v)}{g(G^{-1}(v))},\ \ \
x\in\mathbb{R}^N. \end{eqnarray}
If $v\in H^1(\mathbb{R}^N)\cap L^\infty(\mathbb{R}^N)$ is a solution of $(\ref{S2.00})$, then  it satisfies
$$\int_{\mathbb{R}^N}\left[ \nabla v \nabla\phi +\kappa V_\varepsilon(x)\frac{G^{-1}(v)}{g(G^{-1}(v))}\phi-K_\varepsilon(x)\frac{|G^{-1}(v)|^{p-2}G^{-1}(v)}{g(G^{-1}(v))}\phi \right]dx=0,\quad \forall\phi \in H^1(\mathbb{R}^N).$$
Then, $u=G^{-1}(v)\in H^1(\mathbb{R}^N)\cap L^\infty(\mathbb{R}^N)$.  For any $\varphi\in C_0^\infty(\mathbb{R}^N)$, one has $\varphi g(G^{-1}(v))\in H^1(\mathbb{R}^N)\cap L^\infty(\mathbb{R}^N)$ and satisfies
$$\int_{\mathbb{R}^N}\left[ \nabla u \nabla\varphi +g (u)g' (u)|\nabla u|^2\varphi +\kappa V_\varepsilon(x)u\varphi-K_\varepsilon(x)|w|^{p-2}w\varphi\right] dx=0,$$
which implies that $u$ is a weak solution of (\ref{main}).

\noindent Therefore, in order to find nontrivial solutions to  (\ref{main}), we are renconducted to find nontrivial solutions of  (\ref{S2.00}). Since we are concerned with positive solutions, we actually consider the  following truncated energy functional
\begin{eqnarray*}v\mapsto \frac{1}{2}\int_{\mathbb{R}^N}\left(|\nabla v|^2+\kappa V_\varepsilon(x)|G^{-1}(v)|^2\right)dx-\frac{1}{p}\int_{\mathbb{R}^N}K_\varepsilon(x)|G^{-1}(v^+)|^pdx.
\end{eqnarray*}
However, in order to avoid cumbersome notations, hereafter we write $v$ in place of $v^+$ in the last integral, when this does not yield confusion.

\noindent Set
 \begin{eqnarray}\label{s2.1}\Gamma_{\varepsilon}(v):=P_{\varepsilon} (v)+Q_{\varepsilon}(v),\end{eqnarray}  where
 \begin{eqnarray}\label{2.2}Q_\varepsilon(v)=\Big(\int_{\mathbb{R}^N}\chi_\varepsilon
 v^{\frac{p}{2}}dx-1\Big)_+^{2} \end{eqnarray}
with $\chi_\varepsilon(x)=0$ for $x\in \mathcal{O}_\varepsilon:=\{x\in \mathbb{R}^N:\varepsilon x\in  \mathcal{O}\}$ and $\chi_\varepsilon(x)=\varepsilon^{-\tau}$ for $x\not\in \mathcal{O}_\varepsilon$, where $\tau>0$ has to be determined later on.  By inspection $\Gamma_{\varepsilon}
\in C^1(H^1(\mathbb{R}^N),\mathbb{R})$.
The functional $Q_\varepsilon$ will act as a penalization to force the concentration phenomena
to occur inside $\mathcal{O}$. This type of penalization was introduced in \cite{ByeJea03,ByeJea07}.

\noindent Let $U$ be the unique fast decay positive radial solution  of (\ref{mass1}). Without loss of generality, we may assume $U(0)=\max U(x)$ and that $0\in \mathcal{M}$. Set $
U_{t}(x):=U(\frac{x}{t})$  for $t>0$.
By (\ref{Poh1}), there exists $t_0>1$ such that
\begin{eqnarray}\label{Poho}L_m(U_t)=\left(\frac{t^{N-2}}{2}-\frac{t^N}{2^*}\right)\int_{\mathbb{R}^N}|\nabla U|^2dx<-2,\ \ \text{ for all } t\geq t_0.\end{eqnarray}

\noindent Choose a positive number $\beta <\frac{\text{dist}(\mathcal{M}, \mathbb{R}^N\setminus \mathcal{O})}{100}$ and a cut-off funtion
$\varphi(x)\in C_0^\infty(\mathbb{R}^N,[0,1])$ such that
$\varphi(x)=1$ for $|x|\leq \beta$ and $\varphi(x)=0$ for $|x|\geq
2\beta.$ Set $\varphi_\varepsilon(x):=\varphi(\varepsilon x)$
and then define
$$U_\varepsilon ^{y}(x):=\varphi_\varepsilon \left(x-\frac{y}{\varepsilon}\right)U\left(x-\frac{y}{\varepsilon}\right),\ \  \text{ for each } y\in \mathcal{M}^\beta,$$
where $\mathcal{M}^\beta=\{x\in \mathbb{R}^N: \text{dist}(x,\mathcal{M})<\beta\}$.

\noindent We aim at finding  a solution of (\ref{maineq2} ) near the set
$X_{\varepsilon}:=\{U_\varepsilon ^{y}(x): y\in \mathcal{M}^\beta\}$ for sufficiently small $\varepsilon>0$.
Let
$W_{\varepsilon,t}(x)=\varphi_\varepsilon
(x)U_t(x)$. Note that for fixed $x\in  \mathbb{R}^N$, $W_{\varepsilon,t}(x)\rightarrow 0$ as $t\rightarrow 0$.  So, we  set $ W_{\varepsilon,0}(x)=0.$

\noindent Next we borrow some ideas from  \cite{CWZ20}, however, here the situation is quite different in particular for the decaying behavior of the ground state solution of the limit equation and the different concentrating behavior of the solution.

 \begin{lemma}\label{1}
$ \lim\limits_{\varepsilon\rightarrow 0}\max\limits_{t\in (0,t_0]}|\Gamma_{\varepsilon}(W_{\varepsilon,t})-L_{m}(U_t)|\rightarrow 0.$
\end{lemma}
\begin{proof} Since   $\text{supp}(W_{\varepsilon,t}(x)) \subset \mathcal{O}_\varepsilon$, one has  $Q_\varepsilon(W_{\varepsilon,t}(x))=0$. Thus,
 for $t\in (0,t_0]$, we have
 \begin{eqnarray} \label{lm2.0} \begin{split}
\Gamma_{\varepsilon}(W_{\varepsilon,t})-L_{m}(U_t)=&\frac{1}{2}\int_{\mathbb{R}^N}(|\nabla
W_{\varepsilon,t}|^2-|\nabla U_t|^2)dx+\frac{\kappa}{2}\int_{\mathbb{R}^N}V_\varepsilon|G^{-1}(W_{\varepsilon,t})|^2dx\\
&-\frac{1}{p}\int_{\mathbb{R}^N}(K_\varepsilon|G^{-1}(W_{\varepsilon,t})|^p-m|G^{-1}(U_t)|^p)dx.\end{split}\end{eqnarray}
By  (\ref{estimate}) and the Lebesgue dominated convergence theorem, we get
\begin{eqnarray} \label{r1} \begin{split}
\left|\int_{\mathbb{R}^N}(|\nabla
W_{\varepsilon,t}|^2-|\nabla U_t|^2)dx\right|  \leq& C \varepsilon^2 \int_{\mathbb{R}^N} \left(|\nabla U|^2+U^2\right) dx \\&+Ct_0^{N}\int_{\mathbb{R}^N}|\varphi^2_\varepsilon(tx)-1|(1+|x|)^{2-2N}dx\rightarrow 0\ .
\end{split}
\end{eqnarray}
as $\varepsilon\rightarrow 0$. Clearly, by Lemma \ref{G}$-(iii)$, we obtain
\begin{eqnarray} \label{r2}  \begin{split}
\int_{\mathbb{R}^N}V_\varepsilon|G^{-1}(W_{\varepsilon,t})|^2dx\leq C\int_{\mathbb{R}^N}|W_{\varepsilon,t}|^2dx\leq C\int_{\mathbb{R}^N}U^2dx<+\infty.
\end{split}
\end{eqnarray}

\noindent By the mean value theorem  and dominated convergence again, we obtain
\begin{eqnarray}\label{r5}  \begin{split}\left| \int_{\mathbb{R}^N}(|G^{-1}(U_t)|^p-|G^{-1}(W_{\varepsilon,t})|^p)dx\right|\leq &
p\int_{\mathbb{R}^N}\frac{|G^{-1}(U_t+\theta W_{\varepsilon,t})|^{p-1}}{g(G^{-1}(U_t+\theta W_{\varepsilon,t}))}|U_t-W_{\varepsilon,t}|dx\\
\leq& C\int_{\mathbb{R}^N}(1-\varphi_\varepsilon(t x))(U^2+U^{\frac{p}{2}})dx\\ \rightarrow &0,\ \text{as } \varepsilon\rightarrow 0,\end{split}
\end{eqnarray} where $0<\theta<1.$
Similarly, we have
\begin{eqnarray}\label{r6}  \begin{split}&\left| \int_{\mathbb{R}^N}(K_\varepsilon(x)-m)|G^{-1}(W_{\varepsilon,t})|^pdx\right|\\ &\leq
 t_0^N\int_{\mathbb{R}^N}(K_\varepsilon(tx)-m)\varphi_\varepsilon(tx) U^{\frac{p}{2}}(x)dx\rightarrow &0,\ \text{as } \varepsilon\rightarrow 0.\end{split}
\end{eqnarray}
The desired
conclusion follows from (\ref{lm2.0})--(\ref{r6}).
\end{proof}

\noindent Now, from (\ref{Poho}) and Lemma \ref{1}, there exists $\varepsilon_0>0$, such that
for $\varepsilon\in (0,\varepsilon_0)$,
$$\Gamma_{\varepsilon}(W_{\varepsilon,t_0})\leq L_{m}(U_{t_0})+1<-1.$$
Define the minimax level
$$C_{\varepsilon}=\inf_{\eta_\varepsilon\in\Phi_{\varepsilon}}\max_{s\in[0,1]}\Gamma_{\varepsilon}(\eta_\varepsilon(s)),$$
where $\Phi_{\varepsilon}=\{\eta_\varepsilon\in C([0,1], H^1(\mathbb{R}^N)):\eta_\varepsilon(0)=0,\eta_\varepsilon(1)=W_{\varepsilon,t_0}\}.$
\begin{lemma} \label{CCm}
$\lim\limits_{\varepsilon\rightarrow 0} C_{\varepsilon}=C_m$.
\end{lemma}
\begin{proof}Let $\eta_\varepsilon(s)=W_{\varepsilon,st_0}$, $s\in[0,1]$, such that $\eta_\varepsilon(s)\in \Phi_{\varepsilon}$.
Since $t_0>1$,  from Lemma \ref{1}, we have
\begin{eqnarray*}   \begin{split}\limsup_{\varepsilon\rightarrow 0}C_{\varepsilon}\leq \limsup_{\varepsilon\rightarrow 0}\max_{t\in[0,t_0]}\Gamma_{\varepsilon}(W_{\varepsilon,t})\leq \max_{t\in[0,t_0]}L_{m}(U_t)=C_m.\end{split}
\end{eqnarray*}
It remains to   prove that $\liminf\limits_{\varepsilon\rightarrow 0} C_{\varepsilon}\geq C_m$.  By definition of $C_{\varepsilon}$,
for any $\widetilde{\varepsilon}>0,$  there exists $\tilde{\eta}_\varepsilon\in \Phi_\varepsilon$ such that
\begin{eqnarray}\label{Cm+}\max_{s\in[0,1]}\Gamma_{\varepsilon}(\tilde{\eta}_\varepsilon(s))<C_{\varepsilon}+\widetilde{\varepsilon}.\end{eqnarray}
Since $P_{\varepsilon}(\tilde{\eta}_\varepsilon(0))=0$ and $P_{\varepsilon}(\tilde{\eta}_\varepsilon(1))\leq \Gamma_{\varepsilon}(W_{\varepsilon,t_0})<-1$, there exists $s_0\in (0,1)$ such that
$$P_{\varepsilon}(\tilde{\eta}_\varepsilon(s_0))=-1 \ \ \text{and}\ \ P_{\varepsilon}(\tilde{\eta}_\varepsilon(s))> -1,\ \ s\in[0,s_0).$$
Then, we have
\begin{eqnarray*} \label{3.7}Q_{\varepsilon}(\tilde{\eta}_\varepsilon(s))\leq \Gamma_{\varepsilon}(\tilde{\eta}_\varepsilon(s))+1<C_{\varepsilon}+\widetilde{\varepsilon}+1,\ \ s\in[0,s_0].\end{eqnarray*}
By Lemma \ref{G}$-(v)$, we have
\begin{eqnarray*}\label{lm2.4}  \begin{split} \int_{\mathbb{R}^N\setminus \mathcal{O}_\varepsilon}|G^{-1}(\tilde{\eta}_\varepsilon(s))|^pdx
\leq& \sqrt[4]{2^p}\int_{\mathbb{R}^N\setminus \mathcal{O}_\varepsilon}|\tilde{\eta}_\varepsilon(s)|^{\frac{p}{2}}dx\\ \leq&   \sqrt[4]{2^p}\varepsilon^\tau \left[\sqrt{Q_{\varepsilon}(\tilde{\eta}_\varepsilon(s))} +1\right]  \leq
  \sqrt[4]{2^p}\varepsilon^\tau \left[\sqrt{C_{\varepsilon}+\widetilde{\varepsilon}+1} +1\right],
  \end{split}
\end{eqnarray*}
for $s\in[0,s_0]$. Therefore, the following lower bound holds
\begin{eqnarray} \label{lm3.1-2} \begin{split}
P_{\varepsilon}(\tilde{\eta}_\varepsilon(s))\geq & L_{m}( \tilde{\eta}_\varepsilon(s) )+\frac{1}{p}(m-K_0)\int_{\mathbb{R}^N\setminus \mathcal{O}_\varepsilon}|G^{-1}(\tilde{\eta}_\varepsilon(s))|^pdx\\
 \geq &L_{m}(\tilde{\eta}_\varepsilon(s))+\frac{1}{p}(m-K_0) \sqrt[4]{2^p}\varepsilon^\tau \left[\sqrt{C_{\varepsilon}+\widetilde{\varepsilon}+1} +1\right],\ \ s\in [0,s_0].
\end{split}
\end{eqnarray}
In particular we have
\begin{eqnarray*}  \begin{split}
L_{m}( \tilde{\eta}_\varepsilon(s_0)) \leq\frac{1}{p}(K_0-m)  \sqrt[4]{2^p}\varepsilon^\tau \left[\sqrt{C_{\varepsilon}+\widetilde{\varepsilon}+1} +1\right]-1<0,\ \ \text{for small } \varepsilon>0.
\end{split}
\end{eqnarray*}
Hence,  $\tilde{\eta}_\varepsilon(ts_0)\in \Phi $ and
$\max\limits_{t\in [0,1]}L_m(\tilde{\eta}_\varepsilon(ts_0))\geq C_m.$ So,  by (\ref{Cm+}) and (\ref{lm3.1-2}), we get
\begin{eqnarray*}  \begin{split}
C_{\varepsilon}+\widetilde{\varepsilon}>\max_{s\in[0,s_0]}\Gamma_{\varepsilon}(\tilde{\eta}_\varepsilon(s))\geq C_m+\frac{1}{p}(m-K_0) \sqrt[4]{2^p} \varepsilon^\tau \left[\sqrt{C_{\varepsilon}+\widetilde{\varepsilon}+1} +1\right],
\end{split}
\end{eqnarray*}which yields $\liminf\limits_{\varepsilon\rightarrow 0} C_{\varepsilon}\geq C_m$ since $\widetilde{\varepsilon}$ is arbitrary.
This completes the proof of the lemma.

\end{proof}
\begin{remark} \label{remark1}
Let $D_{\varepsilon}=\max\limits_{s\in[0,1]}\Gamma_{\varepsilon}(W_{\varepsilon,st_0})$. From the proof of  Lemma \ref{CCm} we get $$\lim\limits_{\varepsilon\rightarrow 0}D_{\varepsilon}= C_m\ .$$
\end{remark}

\noindent Next, we consider the space $E_{\varepsilon}^R:=H^1_0(B_{\frac{R}{\varepsilon}}(0))$ endowed with the norm
$$\|v\|_{\varepsilon,R}=\left[\int_{B_{\frac{R}{\varepsilon}}(0)} (|\nabla v|^2+v^2)dx\right]^{\frac{1}{2}} \ .$$
Note that any $v\in E_{\varepsilon}^R$ can be regarded as an element of $H^1(\mathbb{R}^N)$ by defining $v=0$ on $\mathbb{R}^N\setminus B_{\frac{R}{\varepsilon}}(0).$

\noindent Define also the level sets
$$\Gamma_{\varepsilon}^c:=\{u\in E_{\varepsilon}^R: \Gamma_{\varepsilon}(u)\leq c\}$$
and
$$X^d:=\{u\in E_{\varepsilon }^R:\inf\limits_{v\in X}\|u-v\|_{\varepsilon,R}<d\},\quad d>0\ .$$

\noindent In what follows, for small $d>0$, let $v_n\in X_{\varepsilon_n}^d\cap E_{\varepsilon_n}^{R_n}$ with $\varepsilon_n\rightarrow 0$ and $R_n\rightarrow+\infty$ be such that
$$\lim_{n\rightarrow\infty} \Gamma_{\varepsilon_n}(v_n)\leq C_m \ \ \ \mbox{and} \ \lim_{n\rightarrow\infty} \|\Gamma_{\varepsilon_n}'(v_n)\|_{(E_{\varepsilon_n}^{R_n})'}=0.$$

\noindent From the definition of $X_{\varepsilon_n}^d$, we can find   a sequence $\{y_n\}\subset \mathcal{M}^\beta$ such that
\begin{eqnarray}\label{s2.0}\left\|v_n-\varphi_{\varepsilon_n} \left(\cdot-\frac{y_n}{\varepsilon_n}\right)U\left(\cdot-\frac{y_n}{\varepsilon_n}\right)\right\|_{\varepsilon_n,R_n}\leq d.\end{eqnarray}
This implies that $\{v_n\}$ is bounded in $H^1(\mathbb{R}^N)$. Since $\mathcal{M}^\beta$ is compact, we may assume, up to a subsequence, that $y_n\rightarrow y_0  \in \mathcal{M}^\beta$.

\begin{lemma} \label{Lem3.3}
\begin{eqnarray*} \label{nonvanish} \begin{split}
\lim\limits_{n\rightarrow \infty}\sup\limits_{z\in\{z\in \mathbb{R}^N: \frac{1}{2}\beta\leq |\varepsilon_nz-y_n|\leq 3\beta\}}\int_{B_R(z)}v_n^2dx=0,\ \ \ \ \forall R>0.
\end{split}
\end{eqnarray*}
\end{lemma}
\begin{proof}

Suppose by contradiction that there exist $R>0$ and a sequence $\{z_n\}\subset \{z\in \mathbb{R}^N: \frac{1}{2}\beta\leq |\varepsilon_nz-y_n|\leq 3\beta\}$ such that
\begin{eqnarray}  \label{s2.2}
\lim\limits_{n\rightarrow \infty}\int_{B_R(z_n)}v_n^2dx>0.
\end{eqnarray} Assume $\varepsilon_n z_n\rightarrow z_0\in \{z\in \mathbb{R}^N: \frac{1}{2}\beta\leq |z-y_0|\leq 3\beta\}$.
Let $\tilde{v}_n(\cdot):=v_n(\cdot+z_n)$ be such that $\tilde{v}_n\rightharpoonup \tilde{v}$ in $H^1(\mathbb{R}^N)$, $\tilde{v}_n\rightarrow \tilde{v}$ in $L^p_{loc}(\mathbb{R}^N)$, $p\in[2,2^*)$. Then, by (\ref{s2.2}), we get
\begin{eqnarray}  \label{s2.3}
\int_{B_R(0)}|\tilde{v}|^2dx=\lim\limits_{n\rightarrow \infty}\int_{B_R(0)}\tilde{v}_n^2dx>0,
\end{eqnarray}
which yields $\tilde{v}\not\equiv 0.$

 \noindent Let $\phi \in C_0^{\infty}( \mathbb{R}^N)$,  then for large $n$, we have $\phi(\cdot-z_n)\in E_{\varepsilon_n}^{R_n}$. Since $\lim\limits_{n\rightarrow\infty} \|\Gamma_{\varepsilon_n}'(v_n)\|_{(E_{\varepsilon_n}^{R_n})'}=0$, we obtain
\begin{eqnarray} \label{lm3.17} \begin{split}
o_n(1)\|\phi\|_{\varepsilon_n,R_n}=&\int_{\mathbb{R}^N}\left[\nabla \tilde{v}_n\nabla \phi
-K_{\varepsilon_n}(x+z_n)\frac{|G^{-1}(\tilde{v}_n)|^{p-2}G^{-1}(\tilde{v}_n)}
{g(G^{-1}(\tilde{v}_n))}\phi\right]dx\\
&+\kappa_n\int_{\mathbb{R}^N} V_{\varepsilon_n}(x+z_n)\frac{G^{-1}(\tilde{v}_n)}{g
(G^{-1}(\tilde{v}_n))}\phi dx\\
&-p\left(\int_{\mathbb{R}^N}\chi_{\varepsilon_n} |v_n|^{\frac{p}{2}}dx-1\right)_+\int_{\mathbb{R}^N}\chi_{\varepsilon_n}(x+z_n) |\tilde{v}_n|^{\frac{p}{2}-2}\tilde{v}_n\phi dx,
\end{split}
\end{eqnarray} where $\kappa_n=\varepsilon_n^{\frac{2(p-2)\gamma}{4+(p-2)\gamma}}$.

\noindent Clearly,
\begin{eqnarray*} \lim_{n\rightarrow \infty} \int_{\mathbb{R}^N}\chi_{\varepsilon_n}(x+z_n) |\tilde{v}_n|^{\frac{p}{2}-2}\tilde{v}_n\phi dx=0.
\end{eqnarray*} Since $\{v_n\}$ is bounded in $H^1(\mathbb{R}^N)$, we have
\begin{eqnarray} \label{lm3.18}\int_{\mathbb{R}^N} V_{\varepsilon_n}(x+z_n)\frac{G^{-1}(\tilde{v}_n)}{g
(G^{-1}(\tilde{v}_n))}\phi dx<+\infty.
\end{eqnarray}

\noindent By the Lebesgue dominated convergence theorem, we get
\begin{eqnarray} \label{lm3.19} \begin{split}&
\lim_{n\rightarrow \infty}\int_{\mathbb{R}^N}K_{\varepsilon_n}(x+z_n)\frac{|G^{-1}(\tilde{v}_n)|^{p-2}G^{-1}(\tilde{v}_n)}{g
(G^{-1}(\tilde{v}_n))}\phi dx\\&=\int_{\mathbb{R}^N}K(z_0)\frac{|G^{-1}(\tilde{v})|^{p-2}G^{-1}(\tilde{v})}{g
(G^{-1}(\tilde{v}))}\phi dx.
\end{split}
\end{eqnarray}

\noindent Combine (\ref{lm3.17}), (\ref{lm3.18}) and (\ref{lm3.19}),  to have the following
\begin{eqnarray*}  \begin{split}
\int_{\mathbb{R}^N}\left[\nabla \tilde{v} \nabla \phi-K(z_0)\frac{|G^{-1}(\tilde{v})|^{p-2}G^{-1}(\tilde{v})}{g
(G^{-1}(\tilde{v}))}\phi \right]dx=0,\ \ \ \ \forall \phi\in C_0^\infty( \mathbb{R}^N),
\end{split}
\end{eqnarray*}
which implies that $\tilde{v}$ is a positive solution of the following equation
\begin{eqnarray} \label{v1}  \begin{split}
-\Delta \tilde{v}=K(z_0)\frac{|G^{-1}(\tilde{v})|^{p-2}G^{-1}(\tilde{v})}{g
(G^{-1}(\tilde{v}))},\ \ \ \ x\in \mathbb{R}^N.
\end{split}
\end{eqnarray}
Recall that in the right hand side of the above equality, $\tilde{v}$  is actually $\tilde{v}^+$. Thus, by the maximum principle, we get $ \tilde{v}>0$.  Because of $K(z_0)\leq m$, we get $C_{K(z_0)}\geq C_m$.

\noindent Choosing $R>0$ sufficiently large, by Pohozaev's identity we obtain
\begin{eqnarray} \label{C1}  \begin{split}
2\lim\limits_{n\rightarrow \infty}\int_{B_R(z_n)}|\nabla v_n|^2dx\geq \int_{\mathbb{R}^N}|\nabla \tilde{v}|^2dx=NL_{K(z_0)}(\tilde{v})\geq NC_{K(z_0)}\geq NC_m.
\end{split}
\end{eqnarray}
However,  it follows  from (\ref{s2.0}) that
\begin{eqnarray}\label{C2}  \begin{split}
\int_{B_R(z_n)}|\nabla v_n|^2dx\leq& d^2+2\varepsilon_{n}^2\int_{B_R(z_n)} \left| \nabla \varphi_{\varepsilon_n}\left(x-\frac{y_n}{\varepsilon_n}\right)U\left(x-\frac{y_n}{\varepsilon_n}\right)\right|^2 dx\\
&+\int_{B_R(z_n)}\left|\varphi_{\varepsilon_n}\left(x-\frac{y_n}{\varepsilon_n}\right)\nabla U\left(x-\frac{y_n}{\varepsilon_n}\right)\right|^2dx\\
\leq& d+C\varepsilon_n^2+C\int_{B_R(0)}\left(1+\left|x+z_n-\frac{y_n}{\varepsilon_n}\right|\right)^{2-2N}dx.
\end{split}
\end{eqnarray}
Note that $\lim\limits_{n\rightarrow \infty}|z_n-\frac{y_n}{\varepsilon_n}|=+\infty$. Thus, for $n$ large enough,  by (\ref{C1}) and (\ref{C2}), we get
  a contradiction for small $d>0.$  This completes the proof of Lemma \ref{Lem3.3}.

\end{proof}

\noindent Now choose $\eta\in C_0^\infty(\mathbb{R}^N)$ such that $0\leq \eta\leq1$, $\eta(z)=1$ if $z\in \{z\in \mathbb{R}^N: \beta\leq |z|\leq 2\beta \}$,
$\eta(z)=0$ if $z\in \mathbb{R}^N\setminus  \{z\in \mathbb{R}^N: \frac{1}{2}\beta\leq |z|\leq 3\beta \}$. Setting $\eta_n(z)=\eta(\varepsilon_nz-y_n)v_n$, clearly, $\eta_n$ is bounded in $H^1(\mathbb{R}^N)$. Thus, from Lemma \ref{Lem3.3}, we have
\begin{eqnarray} \label{nonvanish1} \begin{split}
\lim\limits_{n\rightarrow \infty}\sup\limits_{z\in \mathbb{R}^N}\int_{\mathbb{R}^N}|\eta_n|^2dx=0.
\end{split}
\end{eqnarray}
 This fact together with Lions' concentration-compactness lemma give  $\eta_n\rightarrow 0$ in $L^{q}(\mathbb{R}^N)$, $q\in (2,2^*)$. So, we obtain
 \begin{eqnarray} \label{nonvanish2} \begin{split}
\lim\limits_{n\rightarrow \infty}\int_{\{x\in \mathbb{R}^N: \beta\leq |\varepsilon_nx-y_n|\leq 2\beta\}}|v_n|^{q}dx
\leq \lim\limits_{n\rightarrow \infty}\int_{\mathbb{R}^N}|\eta_n|^{q}dx
=0.
\end{split}
\end{eqnarray}

\noindent Set $v_{n,1}(\cdot)=\varphi_{\varepsilon_n}(\cdot-\frac{y_n}{\varepsilon_n})v_n(\cdot)$ and $v_{n,2}=v_n-v_{n,1}$. Then, let us prove the following

\begin{lemma} \label{Lem3.41}
$
\Gamma_{\varepsilon_n}(v_{n})\geq \Gamma_{\varepsilon_n}(v_{n,1})+\Gamma_{\varepsilon_n}(v_{n,2})+o_n(1).
$
\end{lemma}
\begin{proof}
 Since $supp(v_{n,1}) \subset \mathcal{O}_\varepsilon$, we have $Q_{\varepsilon_n}(v_{n,1})=0$ and
$Q_{\varepsilon_n}(v_{n,2})=Q_{\varepsilon_n}(v_{n})$.
Therefore,  by Lemma \ref{G}$-(iv)-(v)$ and  $G^{-1}(0)=0$, for large $n$, we deduce that
\begin{eqnarray} \label{S2.3.24} \begin{split}
&\Gamma_{\varepsilon_n}(v_{n,1})+\Gamma_{\varepsilon_n}(v_{n,2})
\\=&\Gamma_{\varepsilon_n}(v_{n})+\int_{\mathbb{R}^N}\varphi_{\varepsilon_n}(x-\frac{y_n}{\varepsilon_n})\left[\varphi_{\varepsilon_n}\left(x-\frac{y_n}{\varepsilon_n}\right)-1\right]|\nabla v_n|^2dx\\&+\frac{\kappa_n}{2} \int_{\mathbb{R}^N}V_{\varepsilon_n}(x)\left[|G^{-1}(v_{n,1})|^2+|G^{-1}(v_{n,2})|^2-|G^{-1}(v_{n})|^2\right]dx\\
&+\frac{1}{p}\int_{\mathbb{R}^N}K_{\varepsilon_n}(x)\left[|G^{-1}(v_n)|^p-|G^{-1}(v_{n,1})|^p-|G^{-1}(v_{n,2})|^p\right]dx+o_n(1)\\
\leq& \Gamma_{\varepsilon_n}(v_{n})+C  \int_{\{x\in \mathbb{R}^N: \beta\leq |\varepsilon_nx-y_n|\leq 2\beta\}}|v_n|^{\frac{p}{2}}dx+o_n(1).\end{split}
\end{eqnarray}
By (\ref{nonvanish2})  and (\ref{S2.3.24}),  we get
the result.
\end{proof}

\noindent In what follows we use the following notation: $\mathcal{A}_0=\{x\in \mathbb{R}^N: |\varepsilon_n x-y_n|\leq \beta\}$, $\mathcal{A}_1=\{x\in \mathbb{R}^N: |\varepsilon_n x-y_n|\geq 2\beta\}$ and
$\mathcal{A}_2=\{x\in \mathbb{R}^N: \beta\leq |\varepsilon_n x-y_n|\leq 2\beta\}$,

\begin{lemma}\label{Lem3.42} $\Gamma_{\varepsilon_n}(v_{n,2})>0.$
\end{lemma}
\begin{proof}
From (\ref{s2.0}), we get
 \begin{eqnarray} \label{vn2} \begin{split}
\int_{\mathcal{A}_1}( |\nabla v_{n,2}|^2+v_{n,2}^2)dx   =\int_{\mathcal{A}_1} (|\nabla v_{n}|^2+v_n^2)dx \leq d^2.
\end{split}
\end{eqnarray}  Similarly,  we get
 \begin{eqnarray} \label{vn2} \begin{split}
\int_{\mathcal{A}_2} (|\nabla v_{n,2}|^2+v_{n,2}^2)dx   \leq C\int_{\mathcal{A}_2} (|\nabla v_{n}|^2+v_{n}^2)dx
   \leq Cd^2+o_n(1).
\end{split}
\end{eqnarray}
For $n$ large enough, we have $\| v_{n,2}\|_{H^1(\mathbb{R}^N)}\leq Cd$  for small $d>0$. On the other hand, by Lemma   \ref{G}-$(i)-(ii)$,  we get $|G^{-1}(v_{n,2})|^p\leq C|v_{n,2}|^{2^*}$.   Hence
\begin{eqnarray} \label{Res2} \begin{split}
\Gamma_{\varepsilon_n}(v_{n,2})\geq& \frac{1}{2}\int_{\mathbb{R}^N}|\nabla v_{n,2}|^2dx-C\int_{\mathbb{R}^N}|v_{n,2}|^{2^*}dx\\
\geq& \frac{1}{2}\int_{\mathbb{R}^N}|\nabla v_{n,2}|^2dx-Cd^{\frac{2^*-2}{2}}\int_{\mathbb{R}^N}|\nabla v_{n,2}|^{2}dx
\geq \frac{1}{4}\int_{\mathbb{R}^N}|\nabla v_{n,2}|^2dx>0.
\end{split}
\end{eqnarray} This concludes the proof of Lemma \ref{Lem3.42}.
\end{proof}

\noindent Denote  the usual norm in $D^{1,2}_0(B_{\frac{R}{\varepsilon}}(0))$ as follows
$$\|v\|^*_{\varepsilon,R}=\Big(\int_{B_{\frac{R}{\varepsilon}}(0)} |\nabla v|^2dx\Big)^{\frac{1}{2}}.$$

\begin{lemma} \label{Lem3.3+}
For small $d>0$, there exist a sequence $\{z_n\}\subset \mathbb{R}^N$ and $y_0\in \mathcal{M}$  with $\varepsilon_n\rightarrow 0$ and $R_n\rightarrow+\infty$ satisfying, up to a subsequence, the following:
$$\lim_{n\rightarrow \infty}|\varepsilon_n z_n-y_0|=0\ \text{ and } \ \ \lim_{n\rightarrow \infty}\|v_n(\cdot)-\varphi_{\varepsilon_n}(\cdot-z_n)U(\cdot-z_n)\|^*_{\varepsilon_n,R_n}=0\ .$$
\end{lemma}
\begin{proof}
Let $w_n(\cdot):=v_{n,1}(\cdot+\frac{y_n}{\varepsilon_n})$. Then,  $\{w_n\}$ is bounded in $H^1(\mathbb{R}^N)$.
 Thus, up to a subsequence if necessary, we may assume $w_n\rightharpoonup w$ in $H^1(\mathbb{R}^N)$, $w_n\rightarrow  w$ in $L_{loc}^q(\mathbb{R}^N)$, $q\in [2,2^*)$, $w_n\rightarrow  w$  a.e. in $\mathbb{R}^N$.
From (\ref{s2.0}),  for given $R>0$, as $n$ is large enough we get
\begin{eqnarray*}  \begin{split}
d^2\geq \int_{\mathcal{A}_0}\left|v_{n,1}-\varphi_{\varepsilon_n}\left(x-\frac{y_n}{\varepsilon_n}\right)U\left(x-\frac{y_n}{\varepsilon_n}\right)\right|^2dx
\geq  \int_{B_R(0)}\left|w_n-\varphi_{\varepsilon_n}U\right|^2dx.
\end{split}
\end{eqnarray*}
Thus, we have
\begin{eqnarray*}  \begin{split}
\int_{B_R(0)}w^2dx= \lim_{n\rightarrow \infty}\int_{B_R(0)}w_n^2dx \geq C-d^2,
\end{split}
\end{eqnarray*}
which yields $w\not=0.$

\noindent Let $\phi\in C_0^\infty( \mathbb{R}^N)$. Note that  $w_n(x)=v_{n,1}(x+\frac{y_n}{\varepsilon_n})=v_n(x+\frac{y_n}{\varepsilon_n})$ for $x\in supp(\phi)$ and large $n$. Moreover, $supp(w_n(x))\subset \{x\in \mathbb{R}^N: |\varepsilon_n x|\leq 2\beta \}\subset\mathcal{O}$.  Thus, from
$\langle \Gamma_{\varepsilon_n}'(v_n), \phi(\cdot-\frac{y_n}{\varepsilon_n})\rangle=o_n(1)\|\phi\|_{\varepsilon_n,R_n}$ and analogously to the proof of (\ref{lm3.18}) and (\ref{lm3.19}),
 $w$ is a positive solution of the following equation
\begin{eqnarray}\label{y0}  \begin{split}
-\Delta w=K(y_0)\frac{|G^{-1}(w)|^{p-2}G^{-1}(w)}{g
(G^{-1}(w))},\ \ \ \ x\in \mathbb{R}^N.
\end{split}
\end{eqnarray}

\noindent \textit{Claim:}
\begin{eqnarray} \label{nonv1} \begin{split}
\lim\limits_{n\rightarrow \infty}\sup\limits_{z\in \mathbb{R}^N}\int_{B_1(z)}|w_n-w|^2dx=0.
\end{split}
\end{eqnarray}
Indeed,
if (\ref{nonv1}) does not occur, then there exists a sequence $\{z_n\}\subset \mathbb{R}^N$ with $|z_n|\rightarrow +\infty$ such that
$$
\lim\limits_{n\rightarrow \infty}\int_{B_1(z_n)}|w_n-w|^2dx>0.
$$
Thus, we have
$$
\lim\limits_{n\rightarrow \infty}\int_{B_1(z_n)}|w|^2dx=0,
\ \ \ \
\lim\limits_{n\rightarrow \infty}\int_{B_1(z_n)}|w_n|^2dx>0.
$$
We have $|\varepsilon_nz_n|\leq \frac{1}{2}\beta$. In fact, if $|\varepsilon_nz_n|> \frac{1}{2} \beta$, by Lemma  \ref{Lem3.3},  we have
\begin{eqnarray*}  \begin{split} \label{s2.6}
0<\lim\limits_{n\rightarrow \infty}\int_{B_1(z_n)}|w_n|^2dx\leq \lim\limits_{n\rightarrow \infty}\sup\limits_{z\in\{z\in \mathbb{R}^N: \frac{1}{2}\beta\leq |\varepsilon_nz-y_n|\leq 3\beta\}}\int_{B_1(z)}|v_n|^2dx=0,
\end{split}
\end{eqnarray*}
 which is impossible.  Thus, up to a subsequence, we may assume $\varepsilon_n z_n\rightarrow z_0\in \{z\in \mathbb{R}^N: |z|\leq \frac{1}{2}\beta \}$.
Suppose $v_{n,1}(\cdot+z_n+\frac{y_n}{\varepsilon_n})\rightharpoonup v_1(\cdot)$ in $H^1(\mathbb{R}^N)$. As in the proof of (\ref{y0}), we have
\begin{eqnarray} \label{v2}  \begin{split}
-\Delta v_1=K(y_0+z_0)\frac{|G^{-1}(v_1)|^{p-2}G^{-1}(v_1)}{g
(G^{-1}(v_1))},\ \ \ \ x\in \mathbb{R}^N.
\end{split}
\end{eqnarray}
By the maximum principle $v_1>0$.

\noindent Thus, for large $R$, we obtain
\begin{eqnarray*}  \begin{split}
\frac{1}{2}NC_m\leq& \int_{B_R(0)}\left|\nabla v_{n,1}\left(x+z_n+\frac{y_n}{\varepsilon_n}\right)\right|^2dx\\ =&\int_{B_R(z_n+\frac{y_n}{\varepsilon_n})}|\nabla v_{n,1}(x)|^2dx\\
\leq& C\varepsilon_n^2+C\int_{B_R(z_n+\frac{y_n}{\varepsilon_n})}|\nabla v_n|^2dx\\
\leq& C\varepsilon_n^2+ Cd+C\int_{B_R(0)}(1+|x+z_n|)^{2-2N}dx.
\end{split}
\end{eqnarray*}
We get  a contradiction  for  large $n$ and small  $d$ since $|z_n|\rightarrow +\infty$.

\noindent Therefore (\ref{nonv1}) holds and the claim is proved.

\noindent Again by Lions' concentration-compactness lemma, we have $w_n\rightarrow w$ in $L^{q}(\mathbb{R}^N)$, $q\in (2,2^*)$.  As a consequence
\begin{eqnarray} \label{Res0} \begin{split}
\lim_{n\rightarrow \infty}\int_{\mathbb{R}^N}K_{\varepsilon_n}\left(x+\frac{y_n}{\varepsilon_n}\right)|G^{-1}(w_n)|^{p}dx=\int_{\mathbb{R}^N}K(y_0)|G^{-1}(w)|^{p}dx\ .
\end{split}
\end{eqnarray}

\noindent By Fatou's lemma, Lemmas \ref{Lem3.41}, \ref{Lem3.42}  and (\ref{Res2}),  up to a subsequence,  we have
\begin{eqnarray}\label{Res3}  \begin{split}
C_m \geq & \lim_{n\rightarrow \infty}\Gamma_{\varepsilon_n}(v_{n})\\
\geq&\lim_{n\rightarrow \infty} \Gamma_{\varepsilon_n}(v_{n,1})+\lim_{n\rightarrow \infty}\Gamma_{\varepsilon_n}(v_{n,2})\\
\geq &\lim_{n\rightarrow \infty}\Gamma_{\varepsilon_n}(v_{n,1}) \\ =& \lim_{n\rightarrow\infty}\frac{1}{2}\int_{\mathbb{R}^N}(|\nabla w_n|^2+\kappa_n V_{\varepsilon_n}\left(x+\frac{y_n}{\varepsilon_n}\right)|G^{-1}(w_n)|^2)dx\\&- \frac{1}{p}\int_{\mathbb{R}^N}K_{\varepsilon_n}\left(x+\frac{y_n}{\varepsilon_n}\right)|G^{-1}(w_n)|^{p}dx\\ \geq & \frac{1}{2}\int_{\mathbb{R}^N} |\nabla w|^2dx- \frac{1}{p}\int_{\mathbb{R}^N}K (y_0)|G^{-1}(w)|^{p}dx
\geq  C_{K(y_0)} \geq C_m.
\end{split}
\end{eqnarray}
Hence, we get $\lim_{n\rightarrow \infty}\Gamma_{\varepsilon_n}(v_{n,1})=C_m.$ Moreover, we get $K(y_0)=m$
and we see that $w$ is a ground state to (\ref{mass1}). Thus, there exists some $z\in \mathbb{R}^N$ such that $w(\cdot+z)=U(\cdot)$.   By (\ref{Res0}) and (\ref{Res3}), we have
\begin{eqnarray*}  \begin{split}
\lim_{n\rightarrow\infty}  \int_{\mathbb{R}^N}|\nabla w_n|^2dx =  \int_{\mathbb{R}^N}|\nabla w|^2dx.
\end{split}
\end{eqnarray*}

\noindent Let $z_n=z+\frac{y_n}{\varepsilon_n}$, then
\begin{eqnarray*}\|v_{n,1}(\cdot)-\varphi_{\varepsilon_n}(\cdot-z_n)U(\cdot-z_n)\|^*_{\varepsilon_n,R_n}\rightarrow 0.
\end{eqnarray*}
Finally, by (\ref{Res2})  and (\ref{Res3}), we have
\begin{eqnarray*} \begin{split}
0=\lim_{n\rightarrow \infty}\Gamma_{\varepsilon_n}(v_{n,2})
\geq \frac{1}{4}\lim_{n\rightarrow \infty}\int_{\mathbb{R}^N}|\nabla v_{n,2}|^2dx,
\end{split}
\end{eqnarray*}
which yields $\lim\limits_{n\rightarrow \infty}\|v_{n,2}\|^*_{\varepsilon_n,R_n}=0$ and the Lemma is proved.
\end{proof}

\noindent Let $d\in (0,d_0)$ such that Lemmas \ref{Lem3.3}, \ref{Lem3.41}, \ref{Lem3.42} and  \ref{Lem3.3+} hold and define
$$\widetilde{X}_\varepsilon^d:=\left\{u\in E_\varepsilon^R:\inf\limits_{v\in X_\varepsilon}\|u-v\|^*_{\varepsilon,R}<d\right\},\quad d>0\ .$$

\begin{lemma} \label{lm6} For any $d\in (0,d_0)$, there exist positive constants $\delta_d$, $R_d$ and $\varepsilon_d$ such that $$\|\Gamma_{\varepsilon}'(v)\|_{(E_{\varepsilon}^R)'}\geq \delta_d$$ for any $v\in E_{\varepsilon}^R\cap \Gamma_{\varepsilon}^{D_{\varepsilon }}\cap (X_\varepsilon^{d_0}\setminus \widetilde{X}_\varepsilon^d)$, $R\geq R_d$ and $\varepsilon\in (0,\varepsilon_d)$.

\end{lemma}
\begin{proof}
By contradiction, we assume that for some $d\in (0,d_0)$, there exists $\varepsilon_n<\frac{1}{n}$, $R_n>n$ and  $v_n\in E_{\varepsilon_n}^{R_n}\cap \Gamma_{\varepsilon_n}^{D_{\varepsilon }}\cap (X_{\varepsilon_n}^{d_0}\setminus \widetilde{X}_{\varepsilon_n}^d)$
such that  $\|\Gamma_{\varepsilon_n}'(v_n)\|_{(E_{\varepsilon_n}^{R_n})'}< \frac{1}{n}$. By Lemma \ref{Lem3.3+}, there exist a sequence $\{z_n\}\subset \mathbb{R}^N$, $y_0\in \mathcal{M}$  satisfying
$$\lim_{n\rightarrow \infty}|\varepsilon_n z_n-y_0|=0,\ \ \ \ \lim_{n\rightarrow \infty}\|v_n-\varphi_{\varepsilon_n}(\cdot-z_n)U(\cdot-z_n)\|^*_{\varepsilon_n,R_n}=0$$ up to a subsequence. Thus, for large $n$,
$\varepsilon_nz_n\in \mathcal{M}^\beta$,  $\varphi_{\varepsilon_n}(\cdot-z_n)U(\cdot-z_n)\in X_{\varepsilon_n}$ and  $v_n\in \widetilde{X}_{\varepsilon_n}^d$, which contradicts the fact $v_n\in X_{\varepsilon_n}^{d_0}\setminus \widetilde{X}_{\varepsilon_n}^d.$
\end{proof}

\begin{lemma}\label{lm7}  For any given $\delta>0$, there exist small  positive constants $\varepsilon_1$ and $d\leq d_0$ such that $\Gamma_{\varepsilon}(v)>C_m-\delta$ for any $v\in X_\varepsilon^{d}$ and $\varepsilon\in (0,\varepsilon_1)$.
\end{lemma}
\begin{proof} For $v\in X_\varepsilon^{d}$, there exists $y\in \mathcal{M}^\beta$   such that
$U_\varepsilon ^{y}(x):=\varphi_\varepsilon (x-\frac{y}{\varepsilon})U(x-\frac{y}{\varepsilon})\in X_\varepsilon$ and
 $\|v-U_\varepsilon ^{y}(x)\|_\varepsilon\leq d$.  Thus, we get
 \begin{eqnarray*}  \begin{split}
\Gamma_{\varepsilon}(U_\varepsilon^y)-C_m\geq &\frac{1}{2}\int_{\mathbb{R}^N}\left[|\nabla (\varphi_\varepsilon U)|^2-|\nabla U|^2\right]dx
+\frac{m}{p}\int_{\mathbb{R}^N}\left(|G^{-1}(U)|^{p}-|G^{-1}(\varphi U)|^{p}\right)dx.
\end{split}
\end{eqnarray*}
Similarly to  the proof of Lemma \ref{1}, for small $\varepsilon>0$, we have
 \begin{eqnarray} \label{eq3.28}  \begin{split}
\Gamma_{\varepsilon}(U_\varepsilon^y)\geq C_m-\frac{\delta}{2}.
\end{split}
\end{eqnarray}
On the other hand, for $v\in X_\varepsilon^{d}$, by choosing $d$ small enough, we have
 \begin{eqnarray} \label{eq3.29} \begin{split} \Gamma_{\varepsilon}(v)-
\Gamma_{\varepsilon}(U_\varepsilon^y)\geq &-\frac{\delta}{2}.
\end{split}
\end{eqnarray}
The result follows from (\ref{eq3.28}) and (\ref{eq3.29}).

\end{proof}

\begin{lemma} \label{L3.10} For sufficiently small $\varepsilon>0$ and large $R>0$, there exists a sequence $\{v_{\varepsilon,n}^R\}\subset  E_{\varepsilon}^R\cap \Gamma_{\varepsilon}^{D_\varepsilon}\cap X_\varepsilon^{d_0}$,  such that $\|\Gamma_{\varepsilon}'(v_{\varepsilon,n}^R)\|_{(E_{\varepsilon}^R)'}\rightarrow 0$ as $n\rightarrow \infty.$
\end{lemma}
\begin{proof}  The proof is similar to   \cite{Gloss0,Gloss01}. For reader's convenience, let us give a detailed proof. By contradiction, for small $\varepsilon>0$   and large $R>0$, there exists $C(\varepsilon,R)>0$ such that
$$\|\Gamma_{\varepsilon}'(v)\|_{(E_{\varepsilon}^R)'}\geq C(\varepsilon,R),\ \ \ v\in E_{\varepsilon}^R\cap \Gamma_{\varepsilon}^{D_\varepsilon}\cap X_\varepsilon^{d_0}.$$
On the other hand, by Lemma \ref{lm6}, there exists $\delta>0$ independent of $\varepsilon\in (0,\varepsilon_0)$ and $R>R_0$ such that $$\|\Gamma_{\varepsilon}'(v)\|_{(E_{\varepsilon}^R)'}\geq \delta,\ \ \ v\in E_{\varepsilon}^R\cap \Gamma_{\varepsilon}^{D_\varepsilon}\cap (X_\varepsilon^{d_0}\setminus \widetilde{X}_\varepsilon^{d_1}).$$
Thus, there exists a pseudo-gradient vector field $\Upsilon_\varepsilon^R$ in a  neighborhood $N_{\varepsilon }^R\subset E_{\varepsilon }^R$ of $E_{\varepsilon }^R\cap \Gamma_{\varepsilon}^{D_\varepsilon}\cap  X_\varepsilon^{d_0} $.
Let $\widetilde{N}_{\varepsilon}^R\subset N_{\varepsilon }^R$ such that $$\|\Gamma_{\varepsilon}'(v)\|_{(E_{\varepsilon }^R)'}\geq \frac{1}{2} C(\varepsilon,R),\ \ \ v\in \widetilde{N}_{\varepsilon }^R.$$
We choose two positive Lipschitz continuous functions $\zeta_\varepsilon^R$ and $\xi$ satisfying
$$\zeta_\varepsilon^R(v)=1, \ \ \mbox{ if } v\in  E_{\varepsilon }^R\cap \Gamma_{\varepsilon}^{D_\varepsilon}\cap  X_\varepsilon^{d_0},\ \ \mbox{ and } \ \ \zeta_\varepsilon^R(v)=0,\ \ \mbox{ if }v \in E_{\varepsilon }^R\setminus \widetilde{N}_{\varepsilon}^R, \ \ 0\leq \zeta_\varepsilon^R\leq 1$$
and
$$\xi\leq 1, \ \ \xi(a)=1, \ \ \mbox{ if }\ \ |a-C_m|\leq \frac{1}{2}C_m, \ \ \mbox{ and } \ \ \xi(a)=0,\ \ \mbox{ if } \ \ |a-C_m|\geq C_m.$$
Define
\begin{eqnarray}\label{4.5}\Psi_\varepsilon^R= \begin{cases}
-\zeta_\varepsilon^R(v)\xi(\Gamma_{\varepsilon}(v))\Upsilon_\varepsilon^R, \quad &\mbox{ if } v\in N_\varepsilon^R,\\
0, \quad &\mbox{ if } v\not\in
 E_{\varepsilon }\setminus N_\varepsilon^R.\end{cases}\end{eqnarray}
Then the initial value problem
\begin{eqnarray}\label{4.6}  \begin{cases}
\frac{d}{dt}F_\varepsilon^R(v,t)=\Psi_\varepsilon^R(F_\varepsilon^R(v,t)),\\
F_\varepsilon^R(v,0)=v\end{cases}\end{eqnarray}
 yields a unique global solution $F_\varepsilon^R:E_{\varepsilon }\times[0,+\infty) \rightarrow E_{\varepsilon }^R$. For the properties of $F_\varepsilon^R$, we refer to e.g. \cite{Gloss01,Str90}.
 Let $\eta_\varepsilon(s)=W_{\varepsilon,st_0}=\varphi_\varepsilon U_{st_0}$, $s\in[0,1]$ as before. Then, for the small $d_1>0$, there exists some $\mu>0$ such that if $|st_0-1|\leq \mu$, then
 $$\|\eta_\varepsilon(s)-\varphi_\varepsilon U\|=\|\varphi_\varepsilon (U_{st_0}-U) \|\leq C\| U_{st_0}-U \|\leq d_1,$$ which implies $\eta_\varepsilon(s)\in  X_\varepsilon^{d_1} \subset \widetilde{X}_\varepsilon^{d_1}$ since $0\in \mathcal{M}$.
On the other hand, if $|st_0-1|\geq \mu$, since $t=1$ is the unique maximum point of $L_m(U_t)$ and $\max\limits_{t\geq 0}L_m(U_t)=L_m(U_1)=C_m$, there exists $\rho>0$ such that
  $L_m(U_{st_0})<C_m-2\rho$ for $|st_0-1|\geq \mu$.
  From Lemma \ref{1}, there exists $\varepsilon_1>0$ such that
  $$\max\limits_{s\in (0,1]}|\Gamma_{\varepsilon}(W_{\varepsilon,st_0})-L_{m}(U_{st_0})|<\rho,\ \ \ \varepsilon\in (0,\varepsilon_1).$$
  So, for $|st_0-1|\geq \mu$, we have
\begin{eqnarray}\label{4.8+}  \begin{split}\Gamma_{\varepsilon}(\eta_\varepsilon(s))=&\Gamma_{\varepsilon}(W_{\varepsilon,st_0})\\ \leq& |\Gamma_{\varepsilon}(W_{\varepsilon,st_0})-L_{m}(U_{st_0})|+L_{m}(U_{st_0})<C_m-\rho,\ \ \ \varepsilon\in (0,\varepsilon_1). \end{split}\end{eqnarray}

\noindent Define $\eta_\varepsilon^R(s,t):
 =F_\varepsilon^R(\eta_\varepsilon(s), t)$, $(s,t)\in[0,1]\times [0,+\infty)$.
Since $\Gamma_{\varepsilon}(\eta_\varepsilon(0)),$ $\Gamma_{\varepsilon}(\eta_\varepsilon(1))\not\in (0,2C_m)$, we get $\eta_\varepsilon^R(s,t)\in \Phi_\varepsilon$  for any $t>0$.

\noindent If $|st_0-1|\geq \mu$, by (\ref{4.8+}),  we have
 $$\Gamma_{\varepsilon}(\eta_\varepsilon^R(s,t))\leq \Gamma_{\varepsilon}(\eta_\varepsilon(s))<C_m-\rho,$$ which is impossible by Lemma \ref{CCm}.

\noindent If $|st_0-1|\leq \mu$, we get $\eta_\varepsilon(s)\in \widetilde{X}_\varepsilon^{d_1}$.  In this case one of the following alternatives holds:
\begin{itemize}
 \item[(a)] $\eta_\varepsilon^R(s,t)\in X_\varepsilon^{d_0}$ for all $t>0$;
 \item[(b)] there exists some $t_s>0$ such that $\eta_\varepsilon^R(s,t_s)\not\in X_\varepsilon^{d_0}$.
 \end{itemize}
 If (a) holds, we have
 \begin{eqnarray*}  \begin{split}
 \Gamma_{\varepsilon}(\eta_\varepsilon^R(s,t))= &\Gamma_{\varepsilon}(\eta_\varepsilon(s))+\int_0^t
 \frac{d}{d\tau}\Gamma_{\varepsilon}(\eta_\varepsilon^R(s,\tau))d\tau\\
 \leq &D_{\varepsilon}-\min\{\delta^2,C(\varepsilon,R)^2\}t.
\end{split}
\end{eqnarray*}
Thus, $\lim\limits_{t\rightarrow +\infty}\Gamma_{\varepsilon}(\eta_\varepsilon^R(s,t))=-\infty$, which contradicts to Lemma \ref{lm7}.
So, we have that (b) holds. For any  fixed $s$ with $|st_0-1|\leq \mu$, we find    $t^1_s,$ $ t_s^2>0$ such that $\eta_\varepsilon^R(s,t) \in X_\varepsilon^{d_0}\setminus \widetilde{X}_\varepsilon^{d_1}$ for $t\in[t^1_s,t_s^2]\subset (0,t_{s})$ for  $|t_s^1-t_s^2|>\sigma$ for some $\sigma>0$ dependent of $d_0$ and $d_1$.
Thus, by Remark \ref{remark1}, we get
\begin{eqnarray*}  \begin{split}
 \Gamma_{\varepsilon}(\eta_\varepsilon^R(s,t_{s_0}))\leq  &\Gamma_{\varepsilon}(\eta_\varepsilon(s))+\int_{t^1_s}^{t^2_s}
 \frac{d}{dt}\Gamma_{\varepsilon}(\eta_\varepsilon^R(s,\tau))d\tau\\
 \leq &D_{\varepsilon}-\delta^2(t_\varepsilon^2-t_\varepsilon^1)\\
< &C_m-\frac{1}{2}\delta^2\sigma,\ \ \ t\in [t_s^1,t_s^2],\ \ \ \mbox{ if }\ \ |st_0-1|\leq \mu.
\end{split}
\end{eqnarray*}
Therefore, since $[0,1]$ is compact, by the covering theorem, for all $s\in [0,1]$ with $|st_0-1|\leq \mu$,  we can find $t_\varepsilon^R$ such that $$\Gamma_{\varepsilon}(\eta_\varepsilon^R(s,t_\varepsilon^R))<C_m-\frac{1}{2}\delta^2\sigma,$$
 which is a contradiction to Lemma \ref{CCm} since $\eta_\varepsilon^R(s,t_\varepsilon^R)\in \Phi_\varepsilon.$

\end{proof}

\begin{lemma} \label{existence}
For sufficiently small $\varepsilon>0$, there exists a critical point $v_\varepsilon\in  X_{\varepsilon}^{d_0}\cap \Gamma_{\varepsilon}^{D_\varepsilon}$ of $\Gamma_{\varepsilon}$.
\end{lemma}
\begin{proof}
By Lemma \ref{L3.10},  there exists $\varepsilon_0$ and $R_0>0$ such that   there exists a sequence $\{v_{\varepsilon,n}^R\}\subset  E_{\varepsilon}^R\cap \Gamma_{\varepsilon}^{D_{\varepsilon}}\cap X_\varepsilon^{d_0}$,  such that $\|\Gamma_{\varepsilon}'(v_{\varepsilon,n}^R)\|_{(E_{\varepsilon}^R)'}\rightarrow 0$ as $n\rightarrow \infty$ for $ \varepsilon\in (0,\varepsilon_0)$ and $ R\in (R_0,+\infty).$ Clearly, $\{v_{\varepsilon,n}^R\}$ is bounded in $H_0^1(B_{\frac{R}{\varepsilon}}(0))$ since $v_{\varepsilon,n}^R
\in X_\varepsilon^{d_0}$.  Up to a subsequence if necessary, we may assume $v_{\varepsilon,n}^R\rightharpoonup v_{\varepsilon}^R$ in $H_0^1(B_{\frac{R}{\varepsilon}}(0))$, $v_{\varepsilon,n}^R\rightarrow v_{\varepsilon}^R$ in $L^p(B_{\frac{R}{\varepsilon}}(0))$,  $p\in [1,2^*)$ and $v_{\varepsilon,n}^R\rightarrow  v_{\varepsilon}^R$ a.e. in $\mathbb{R}^N$. Thus $v_{\varepsilon}^R$ is a solution of
\begin{eqnarray}  \begin{split}\label{s4.8}-\Delta v=&\frac{K_\varepsilon(x)|G^{-1} (v)|^{p-2}G^{-1} (v)-\kappa V_\varepsilon(  x)G^{-1} (v)}{g (G^{-1} (v))}\\&-p\left(\int_{B_{\frac{R}{\varepsilon}(0)}}\chi_\varepsilon |v|^{\frac{p}{2}}dx-1\right)_+\chi_\varepsilon |v|^{\frac{p}{2}-2}v,\ \ \ \ \ x\in B_{\frac{R}{\varepsilon}}(0).\end{split}
\end{eqnarray}
From (\ref{s4.8}) we have $v_{\varepsilon,n}^R\rightarrow v_{\varepsilon}^R$ in $H_0^1(B_{\frac{R}{\varepsilon}}(0))$ and $v_{\varepsilon}^R\in X_\varepsilon^{d_0}\cap \Gamma_{\varepsilon}^{D_{\varepsilon}}$. By the maximum principle, $v_{\varepsilon}^R>0$.
Note that  any positive solutions of (\ref{s4.8}) satisfies
$$-\Delta v\leq Cv^{\frac{p}{2}-1}, \ \ x\in B_{\frac{R}{\varepsilon}}(0),$$ where $C>0$ is independent of $\varepsilon$ and $R$. In particular,
$-\Delta v_\varepsilon^R\leq C(v^R_\varepsilon)^{\frac{p}{2}-1}, $ $ x\in B_{\frac{R}{\varepsilon}}(0).$  By applying standard Moser's iteration (see \cite{Gil89}),   $\{v_\varepsilon^R\}$ is bounded in
$L^q_{loc}(\mathbb{R}^N)$ uniformly on $R\geq R_0$ and $\varepsilon\in (0,\varepsilon_0)$ for any $q<\infty$. Moreover, for any $y\in \mathbb{R}^N$, we have
\begin{eqnarray}  \begin{split}\label{e1}\|v_\varepsilon^R\|_{L^q(B_3(y))}\leq C\|v_\varepsilon^R\|_{L^{\frac{p}{2}}(B_4(y))}.\end{split}
\end{eqnarray}
By Theorem 8.17 \cite{Gil89} and (\ref{e1}), we have
\begin{eqnarray}  \begin{split}\label{e2}\sup\limits_{B_1(y)}v_\varepsilon^R\leq C(\|v_\varepsilon^R\|_{L^\frac{p}{2}(B_2(y))}+\|(v_\varepsilon^R)^{\frac{p}{2}-1}\|_{L^q(B_3(y))})\leq C\|v_\varepsilon^R\|_{L^\frac{p}{2}(B_4(y))}.\end{split}
\end{eqnarray}
 In particular this implies that $v_\varepsilon^R$ stays bounded in
$L^\infty(\mathbb{R}^N)$.
Since $\|v_{\varepsilon}^R\|_{\varepsilon}$ and $\{\Gamma_{\varepsilon}(v_{\varepsilon}^R)\}$
are bounded, we get $\{Q_\varepsilon (v_{\varepsilon}^R)\}$ is uniformly bounded on $R\geq R_0$ and $\varepsilon\in (0,\varepsilon_0)$. So, we have $$\int_{\mathbb{R}^N\setminus B_{\frac{R_0}{\varepsilon}}(0)}|v_\varepsilon^R|^{\frac{p}{2}}dx\leq
\int_{\mathbb{R}^N\setminus \mathcal{O}_\varepsilon}|v_\varepsilon^R|^{\frac{p}{2}}dx=\varepsilon^\tau \int_{\mathbb{R}^N}\chi_\varepsilon|v_\varepsilon^R|^{\frac{p}{2}}dx\leq \varepsilon^\tau  C$$ for any $R\geq R_0$ and $\varepsilon\in (0,\varepsilon_0)$.
Thus, for  $|x|\geq \frac{R_0}{\varepsilon}+4$ and $R\geq R_0$, we have
$(v_{\varepsilon}^R)
^{\frac{p}{2}-1}\leq \varepsilon^\tau  C v_{\varepsilon}^R.$ By the comparison principle,
similarly to the proof of Proposition 3 in \cite{ByeJea03}, we get
  \begin{eqnarray} \label{limitA} \begin{split}
\lim_{A\rightarrow +\infty}\int_{\mathbb{R}^N\setminus B_A(0)}\left[|\nabla v_\varepsilon^R|^2+\kappa V_\varepsilon(x)|G^{-1}(v_{\varepsilon}^R)|^2\right]dx=0
\end{split}
\end{eqnarray}
uniformly on $R\geq R_0$. Let $v_k=v_\varepsilon^{R_k}$ and $R_k\rightarrow +\infty$ as $k\rightarrow\infty.$ Then $\{v_k\}$ is bounded in $H^1(\mathbb{R}^N)$ and we may assume $v_k\rightharpoonup v_\varepsilon$ in $H^1(\mathbb{R}^N)$,
$v_k\rightarrow v_\varepsilon$ a.e. in $\mathbb{R}^N$.  Since $v_k$ satisfies (\ref{s4.8}) and using (\ref{limitA}), we get
$\|v_k-v_\varepsilon\|_\varepsilon\rightarrow 0$ as $k\rightarrow\infty.$ Thus, $v_\varepsilon \in X_\varepsilon^{d_0}\cap \Gamma_{\varepsilon}^{D_{\varepsilon}}$ and $\Gamma_{\varepsilon}'(v_\varepsilon)=0.$

\end{proof}

\subsection{Proof of  Theorems \ref{Th1} and \ref{Th2}. }  By   Lemma \ref{existence}, for small $\varepsilon>0$,  there exists a  positive solution $v_\varepsilon$ to the following equation
\begin{align}\label{S4.1} -
\Delta v+\kappa V_\varepsilon(x)\frac{G^{-1}(v)}{g(G^{-1}(v))}
=K_{\varepsilon}(x)\frac{|G^{-1}(v)|^{p-2}G^{-1}(v)}{g (G^{-1} (v))}-p\left(\int_{\mathbb{R}^N}\chi_\varepsilon v^{\frac{p}{2}}dx-1\right)_+\chi_\varepsilon v^{\frac{p}{2}-1}. \end{align}
 Since $v_\varepsilon\in X_\varepsilon^{d_0}$,  by Moser's iteration \cite{Gil89}, $\{v_\varepsilon\}$ is uniformly bounded in $L^\infty(\mathbb{R}^N)$ for small  $\varepsilon>0$.

\noindent By Lemma \ref{Lem3.3}, for small $d>0$, there exist a sequence $\{z_\varepsilon\}\subset \mathbb{R}^N$ and $y_0\in \mathcal{M}$  satisfying
$$\lim_{\varepsilon \rightarrow 0}|\varepsilon z_\varepsilon-y_0|=0,\ \ \ \ \lim_{\varepsilon \rightarrow 0}\|v_{\varepsilon}(\cdot)-\varphi_{\varepsilon}(\cdot-z_\varepsilon)U(\cdot-z_\varepsilon)\|^*_{\varepsilon}=0$$ up to a subsequence. Then,
\begin{align}\label{SS0}  \lim_{n\rightarrow \infty}\|v_{\varepsilon}(\cdot+z_\varepsilon)-U(\cdot)\|^* =0.\end{align}
From (\ref{S4.1}), there exist some $C_1>0$ and $C_2>0$
\begin{align}\label{S2.0+} -
\Delta v_{\varepsilon}+  C_1 \kappa v_{\varepsilon}
\leq C_2v_{\varepsilon}^{2^*-1} . \end{align}
 Thus, for given $\sigma>0$, there exist $R>0$ and $\varepsilon_0>0$ such that
\begin{align}\label{ss1}\sup\limits_{\varepsilon\in (0,\varepsilon_0)}  \kappa\int_{\mathbb{R}^N\setminus B_R(0)}v_{\varepsilon}^2(x+z_{\varepsilon})dx\leq \sigma.\end{align}
Setting $w_\varepsilon(\cdot)=v_\varepsilon(\cdot+z_\varepsilon)$, we have $-\Delta w_\varepsilon\leq Cw_\varepsilon.$
Hence, from Theorem 8.17 in \cite{Gil89}, there exists a constant $C_0=C_0(N,C)$ such that
$$\sup_{B_1(y)}w_\varepsilon\leq C_0\|w_\varepsilon\|_{L^{2}(B_2(y))},\ \ \text{ for all } y\in \mathbb{R}^N.$$ In view of (\ref{ss1}), we conclude that $  w_\varepsilon(x)\rightarrow 0$ as $|x|\rightarrow \infty.$ Let $y_\varepsilon$ be a maximum point of $   w_\varepsilon(x)$, then $\{y_\varepsilon\}$ is bounded.
Otherwise, $\|w_\varepsilon(y_\varepsilon)\|_{L^\infty}\rightarrow 0$ as $\varepsilon\rightarrow 0$, which would contradict  (\ref{SS0}).

\noindent Now, fix $\varepsilon>0$ sufficiently small. By Lemma \ref{G}$-(i)$, we choose $R_0>0$ such that
$$\frac{G^{-1}(v)}{g(G^{-1}(v))}\geq \frac{1}{2}v,\ \ \ \ \text{ for }|x|\geq R_0 $$
Thus, from (\ref{S4.1}), we have
$$-\Delta v+\frac{1}{2}\kappa V_0 v\leq K_0 v^{p-1},\ \ \text{ for }|x|\geq R_0.$$
Let $\phi(y)=\kappa^{\frac{1}{2-p}} v(\frac{y}{\sqrt{\kappa}})$, $|y|\geq \sqrt{\kappa}R_0$, then
$$-\Delta \phi+\frac{1}{2}  V_0 \phi\leq K_0 \phi^{p-1}.$$
Since $\frac{1}{\sqrt{\kappa}}|y|\rightarrow +\infty$, we get $\phi(y)\rightarrow 0,$ thus, there exists $R_1>0$, such that
$$K_0|w|^{p-2}\leq \frac{1}{4}V_0.$$
Thus, for $|y|\geq \sqrt{\kappa}R_1$, we have
$$-\Delta \phi+\frac{1}{4}  V_0 \phi\leq 0.$$

\noindent Now define the function $$\psi(x)=M\exp(-\xi|y|),$$
where $\xi$ and $M$ are such that $4\xi^2<V_0$ and for all $|y|=\sqrt{\kappa}R_1$,
$$M\exp(-\xi \sqrt{\kappa}R_1)>\phi(y)$$
It is straightforward to check that for all $x\not=0$,
$$\Delta \psi\leq \xi^2 \psi \leq \frac{1}{4}V_0 \psi.$$

\noindent Thus
$$-\Delta(\psi-\phi)+\frac{1}{4}V_0(\psi-\phi)\geq 0, \ \ \text{ for }|y|\geq \sqrt{\kappa}R_1.$$
By the maximum principle, we have that
$$\phi(y)\leq M\exp(-\xi|y|),\ \ \ \text{ for }|y|\geq \sqrt{\kappa}R_1,$$
which yields that
\begin{align}\label{ss11}w_\varepsilon(x)\leq M\kappa^{\frac{1}{p-2}}\exp(-\xi \sqrt{\kappa}|x|), \ \ \ \text{ for }|x|\geq R_1.\end{align}

\noindent and in turn
\begin{align}\label{ss12}w_\varepsilon(x)\leq C\exp(-\xi \sqrt{\kappa}|x|), \ \ \ \text{ for } x\in \mathbb{R}^N.\end{align}

Setting $x_\varepsilon=y_\varepsilon+z_\varepsilon$, then $x_\varepsilon$  is a maximum point of $   v_\varepsilon(x)$ and
$$v_\varepsilon(x)=w_\varepsilon(x-z_\varepsilon)\leq C\exp(-\xi \sqrt{\kappa}|x-x_\varepsilon|), \ \ \ \text{ for }x\in \mathbb{R}^N.$$

As a consequence we have that
 \begin{eqnarray*} \begin{split}\int_{\mathbb{R}^N}\chi_\varepsilon
 v_\varepsilon^{\frac{p}{2}}dx=&\varepsilon^{-\tau}\int_{\mathbb{R}^N\setminus \mathcal{O}_\varepsilon}
 v_\varepsilon^\frac{p}{2}dx\\ \leq& C  \varepsilon^{-\tau}\int_{\mathbb{R}^N\setminus \mathcal{O}_\varepsilon}
\exp(-c \sqrt{\kappa}|x-x_\varepsilon|)dx\\
\leq& C  \varepsilon^{-\tau-N}\int_{\mathbb{R}^N\setminus \mathcal{O}}
\exp(-c \sqrt{\kappa}|x/{\varepsilon}-x_\varepsilon|)dx\\
\rightarrow &0, \text{ as } \varepsilon\rightarrow 0.\\
 \end{split}\end{eqnarray*}
(Notice that here we use that fact $\sqrt{\kappa}/\varepsilon=1/\hbar\rightarrow 0$ as $\varepsilon\rightarrow 0$). Thus, $Q_{\varepsilon}(v_\varepsilon)=0$ for small $\varepsilon>0 $ and $v_\varepsilon$ is a positive critical point of $P_{\varepsilon}$. So, $u_{\varepsilon}=G^{-1}(v_\varepsilon)$ is a positive solution of (\ref{maineq2}).
 Furthermore,
$$\|u_{\varepsilon}(\cdot+x_{\varepsilon})-G^{-1}(U(\cdot+y_0))\|_{D^{1,2}(\mathbb{R}^N)} \rightarrow 0.$$

\medskip
\noindent Finally, let us prove the last part of Theorem \ref{Th2} which concerns the critical case $p=\frac{2N}{N-2}$. Set
$$G_\epsilon^{-1}(t)=\int_0^tg_\epsilon(s) ds,$$
where $g_\epsilon(s)=\sqrt{1+2\zeta s^2}$.

\noindent Then, equation (\ref{maineq5}) turns into the following equation
\begin{eqnarray}\label{S2.0} -
\Delta v+\kappa V_\varepsilon(x)\frac{G_\epsilon^{-1}(v)}{g_\epsilon(G_\epsilon^{-1}(v))}
=K_\varepsilon(x)\frac{|G_\epsilon^{-1}(v)|^{p-2}G_\epsilon^{-1}(v)}{g_\epsilon(G_\varepsilon^{-1}(v))},\ \ \
x\in\mathbb{R}^N. \end{eqnarray}

\noindent Here we just stress the differences with respect to the previous case. First, the unique fast decay positive radial solution  of (\ref{mass1}) should be replaced  by the unique positive radial solution  (ground state) of (\ref{limit2}). Besides, $$L_m(v)=\frac{1}{2}\int_{\mathbb{R}^N}|\nabla v|^2dx-\frac{m}{p}\int_{\mathbb{R}^N}|v|^{p}dx.$$

\noindent Since $g_\epsilon(t)\rightarrow 1$ and $G_\epsilon^{-1}(t)\rightarrow t$ for any $t\in \mathbb{R}$ as $\epsilon\rightarrow 0$, (\ref{lm3.19}) turns into
\begin{eqnarray} \label{lm4.19} \begin{split}
\lim_{n\rightarrow \infty}\int_{\mathbb{R}^N}K_{\varepsilon_n}(x+z_n)\frac{|G_\epsilon^{-1}(\tilde{v}_n)|^{\frac{4}{N-2}}G_\epsilon^{-1}(\tilde{v}_n)}{g_\epsilon
(G_\epsilon^{-1}(\tilde{v}_n))}\phi dx=\int_{\mathbb{R}^N}K(z_0)|\tilde{v}|^{\frac{4}{N-2}}\tilde{v}\phi dx.
\end{split}
\end{eqnarray}
So, (\ref{v1}) turns into
\begin{eqnarray} \label{v41}  \begin{split}
-\Delta \tilde{v}=K(z_0)|\tilde{v}|^{\frac{4}{N-2}}\tilde{v},\ \ \ \ x\in \mathbb{R}^N.
\end{split}
\end{eqnarray}

\noindent Similarly, (\ref{v2}) becomes
\begin{eqnarray} \label{V42}  \begin{split}
-\Delta v_1=K(y_0+z_0)|v_1|^{\frac{4}{N-2}}v_1,\ \ \ \ x\in \mathbb{R}^N.
\end{split}
\end{eqnarray}
The rest can be discussed in a similar fashion and following the analysis carried out in the previous case; the proofs are thus complete.

\begin{remark} (\ref{ss11}) is not true for all $x\in \mathbb{R}^N$. In fact,
from  (\ref{S4.1}), we deduce that
 \begin{eqnarray*} \begin{split}\int_{\mathbb{R}^N} |\nabla w_\varepsilon|^2dx<
&\int_{\mathbb{R}^N}\left[1+\frac{g'(G^{-1} (w_\varepsilon))G^{-1} (w_\varepsilon)}{g (G^{-1} (w_\varepsilon))}\right]|\nabla w_\varepsilon|^2dx\\& +\kappa \int_{\mathbb{R}^N}V_\varepsilon(x+z_\varepsilon)|G^{-1}(w_\varepsilon)|^2dx
\\ \leq& \int_{\mathbb{R}^N} K_{\varepsilon}(x+z_\varepsilon) |G^{-1}(w_\varepsilon)|^{p}dx\leq C \|w_\varepsilon\|_{L^\infty(\mathbb{R}^N)}^{p-2^*} \left[\int_{\mathbb{R}^N} |\nabla w_\varepsilon|^2dx\right]^{\frac{2^*}{2}}.  \end{split}\end{eqnarray*}
\end{remark}
By (\ref{SS0}), we get $\|w_\varepsilon\|_{L^\infty(\mathbb{R}^N)}\geq C,$ where $C$ is independent on $\kappa$.

\end{document}